\documentclass[twocolumn]{autart}    
\usepackage{subfig}
\usepackage{graphicx} 
\usepackage{float}
\usepackage{caption}
\usepackage{tikz}
\usepackage{amssymb}  
\usepackage{amsfonts}
\usepackage{amsmath}
\usepackage{mathrsfs} 
\usepackage{amssymb}
\usepackage{mathtools}
\usepackage{commath} 
\usepackage{pgfplots}
\usepackage{siunitx}
\usepackage{mathrsfs} 
\usepackage{xcolor}
\usepackage[noend]{algpseudocode}

\begin{document}
\newtheorem{prop1}{Proposition}
\newtheorem{remark}{Remark}
\newtheorem{lemma}{Lemma}
\newtheorem{problem}{Problem}
\newtheorem{definition}{Definition}
\newcommand{\dd}{\mathop{}\!\mathrm{d}}
\newtheorem{assumption}{Assumption}
\begin{frontmatter}
\title{Reduced order in domain control of distributed parameter port-Hamiltonian systems via energy shaping \thanksref{footnoteinfo}} 
\thanks[footnoteinfo]
{
	This work has been supported by the EIPHI Graduate school (contract ``ANR-17-EURE-0002") and the ANR Project IMPACTS. The fourth author acknowledges Chilean FONDECYT 1191544 and CONICYT BASAL FB0008 projects. Corresponding author Yongxin Wu.
}
\author[Paestum]{Ning Liu}\ead{ning.liu@femto-st.fr},  
\author[Paestum]{Yongxin Wu}\ead{yongxin.wu@femto-st.fr},
\author[Paestum]{Yann Le Gorrec}\ead{yann.le.gorrec@ens2m.fr},
\author[Rome]{Laurent Lefevre}\ead{laurent.lefevre@lcis.grenoble-inp.fr},               
\author[Baiae]{Hector Ramirez}\ead{hector.ramireze@usm.cl}  
\address[Paestum]{Institute FEMTO-ST CNRS UMR 6174, Universit\'e Bourgogne Franche-Comt\'e, 26 chemin de l'\'epitaphe, F-25030, Besan\c{c}on, France}                                              
\address[Rome]{Universit\'e Grenoble Alpes, LCIS, F-26902, Valence, France}            
\address[Baiae]{Universidad Tecnica Federico Santa Maria,Valparaiso, Chile}        
\begin{abstract}                        
	An in-domain finite dimensional controller for a class of distributed parameter systems on a one-dimensional spatial domain formulated under the port-Hamiltonian framework is presented. Based on  \cite{tranchantcdc} where positive feedback and a late lumping approach are used, we extend the Control by Interconnection method and propose a new energy shaping methodology with an early lumping approach on the distributed spatial domain of the system. Our two main control objectives are to stabilize the closed loop system, as well as to improve the closed loop dynamic performances. With the early lumping approach, we investigate two cases of the controller design, the ideal case where each distributed controller acts independently on the spatial domain (fully-actuated), and the more realistic case where the controller is piecewise constant over certain interval (under-actuated). We then analyze the asymptotic stability of the closed loop system when the infinite dimensional plant system is connected with the finite dimensional controller.  Furthermore we provide simulation results comparing the performance of the fully-actuated case and the under-actuated case with an example of an elastic vibrating string.
\end{abstract}
\end{frontmatter}
	
\section{Introduction}
	The control of distributed parameter systems governed by partial differential equations (PDEs) has been the subject of an intensive research activity over the last decades \cite{Curtain1995}. Based on the location of actuators and sensors, control of PDEs is concerned with either boundary or in-domain distributed control. Up to now most of results on control design have been established for boundary controlled systems \cite{krstic2008boundary}, and the study of practical in-domain distributed control is quite limited. However, this latter plays an important role in the field of fluid dynamics,  chemical processes, and recently for flexible structures actuated by soft actuators \cite{Schr_ck_2010,Liu2021TMECH} due to  their recent development. Moreover, it has the potential to be extended to higher dimensional spatial domains.  In this paper, we focus on in-domain distributed control of distributed parameter port-Hamiltonian systems (PHSs), for its advantages in the modeling and control of multi-physical nonlinear systems \cite{duindam2009modeling}. 
	
	The concept of PHSs has firstly been introduced in the 90's in \cite{maschke1992intrinsic} for lumped parameter systems. It is an energy-based representation that expresses the power exchanges within the system and with its environment, using energy and co-energy variables and an intrinsic geometric structure, the Dirac structure. PHSs have later been extended to distributed parameter systems in \cite{maschke2000port}. The geometric Stokes-Dirac structure for linear distributed parameter PHSs defined on a one-dimensional (1D) spatial domain has been investigated in \cite{le2005dirac} using the semigroup theory. In terms of control, the Hamiltonian function is a good Lyapunov candidate function and this class of systems exhibits an intrinsic passivity property, which makes passivity-based control (PBC) techniques suitable. Moreover, controllers designed through PBC have a clear physical interpretation. The PBC has been applied to lumped parameter PHSs and extensively studied in \cite{ortega2001putting,van2000l2} leading to efficient control design techniques such as Control by Interconnection (CbI) and Interconnection and Damping Assignment (IDA)-PBC. \cite{alessandroTAC} has generalized the CbI to distributed parameter PHSs with boundary control. The controller is designed to be a PHS. With the passive interconnection between the plant and the controller, the closed loop system is again a PHS \cite{van2000l2}. Casimir functions are then used to find an invariant relation between the state variables of the plant and those of the controller. By modifying the controller parameters, one can add damping to the closed loop system (Damping Injection) and adjust the shape of its closed loop Hamiltonian modifying both equilibrium and dynamic performances (Energy Shaping). The first result on CbI for distributed parameter PHSs with in-domain distributed control can be found in \cite{tranchantcdc} where a positive feedback, late lumping approach and full actuation are investigated.
	Different from \cite{tranchantcdc}, we consider here the CbI of distributed PHS using a negative feedback, an early-lumping approach and a limited number of actuators. 
	
	The main contribution of this paper are:
	\begin{enumerate}
	\item  The generalization of the in-domain CbI established for the Timoshenko beam in \cite{Liu2021} to a class of 1D linear distributed parameter PHSs which covers different physical applications such as vibrating strings, Timoshenko beams,  Euler-Bernoulli beams, etc.  With the early lumping approach, the finite dimensional controller is designed on the basis of the approximated finite dimensional plant system. 
	\item Using the semigroup theory and passivity properties, we prove the asymptotic stability of the closed loop system when the finite dimensional controller is applied to the infinite dimensional plant getting rid of the well known spillover effect \cite{BALAS1982361}. 
	\item Besides the asymptotic stabilization of the closed loop system, the controller improves the dynamic performances of the system over a given range of frequencies \textit{i.e.} accelerate the system with less oscillation and less overshoot. 
	\end{enumerate}
	Two different cases are investigated for the dynamic performances improvement: the ideal \emph{fully-actuated case} where the control input works independently on each element of the discretized model and the \emph{under-actuated case} where the input acts identically on sets of elements, providing less degrees of freedom. This latter case is closer to the real implementation because the control is usually carried out through patches that act similarly over spatial elements. It is shown how to change the closed-loop energetic properties of the discretized system in a perfect way when the system is fully-actuated and in an \emph{optimal way} when the system is under-actuated.

	The paper is organized as follows: in Section \ref{Secmodelling} the port Hamiltonian formulation of a class of linear distributed parameter systems with two conservation laws and with in-domain control is presented, together with its structure preserving discretization. 
	The aforementioned CbI and energy shaping methods are investigated in Section \ref{CIES} with detailed closed-loop stability analysis. This control strategy takes advantage of both the early lumping approach, leading to a directly implementable controller with guaranteed performances over a given frequency range. The passivity of the system and the Damping Injection guarantee the well-posedness and the asymptotic stability of the closed-loop system. Section \ref{SecSimulation} provides some simulation results with a comparison between the fully- and under-actuated cases using a vibrating string example. Section \ref{SecConclusion} ends up with conclusions and perspectives.
\section{Port-Hamiltonian systems with in-domain control}
	\label{Secmodelling}
	In this paper, we consider partitioned port Hamiltonian systems defined on a one dimensional spatial domain $\zeta\in\left[0,L\right]$ with distributed and boundary control and observation of the form:
	\begin{eqnarray}
\label{twoconservationlaws}
\frac{\partial}{\partial t}\begin{bmatrix}x_1(\zeta,t) \\ x_2(\zeta,t) \end{bmatrix}=\begin{bmatrix}0 & {\mathcal G}\\ -{\mathcal G}^* & -R \end{bmatrix}\begin{bmatrix}\mathcal{L}_1(\zeta) x_1(\zeta,t) \\ \mathcal{L}_2(\zeta) x_2(\zeta,t) \end{bmatrix}\\+\begin{bmatrix}0 \\ B_0\end{bmatrix}u_d(\zeta,t) \\
y_d(\zeta,t)=\begin{bmatrix}0 & B_0^* \end{bmatrix}\begin{bmatrix}\mathcal{L}_1(\zeta) x_1(\zeta,t) \\ \mathcal{L}_2(\zeta) x_2(\zeta,t) \end{bmatrix} \\ \label{output}
u_b= {\mathcal B} \begin{bmatrix}\mathcal{L}_1(\zeta) x_1(\zeta,t) \\ \mathcal{L}_2(\zeta) x_2(\zeta,t) \end{bmatrix}, \;
y_b= {\mathcal C} \begin{bmatrix}\mathcal{L}_1(\zeta) x_1(\zeta,t) \\ \mathcal{L}_2(\zeta) x_2(\zeta,t) \end{bmatrix}
\end{eqnarray}

where $x=[x_1^T, x_2^T]^T \in X \coloneqq L^2(\left[ a,b\right], \mathbb{R}^n) \times L^2(\left[ a,b\right], \mathbb{R}^n)$, 
$\mathcal{L}=\text{diag}(\mathcal{L}_1,\mathcal{L}_2)$ is a bounded and Lipschitz continuous matrix-valued function such that $\mathcal{L}(\zeta)=\mathcal{L}^T(\zeta)$ and $\mathcal{L}(\zeta) \geq \eta$ with $\eta>0$ for all $\zeta \in [a,b]$, $R \in \mathbb{R}^{(n,n)},  R=R^T>0$, ${\mathcal B}(\cdot)$ and ${\mathcal C}(\cdot)$ are some boundary input and boundary output mapping operators that will be defined later. The state space $X$ is endowed with the inner product $\left< x | \tilde{x} \right>_{\mathcal L}=\left< x | {\mathcal L}\tilde{x} \right>$ and norm $\| x \|^2_{\mathcal L}=\left< x | x \right>_{\mathcal L}$ where $\left< \cdot | \cdot  \right>$ denotes the natural $L^2$-inner product. $X\ni x$ is the space of energy variables and ${\mathcal L}x$ denotes the co-energy variable associated to the energy variable $x$. The total energy of the system is given by
\begin{equation}
		H(x_1,x_2)=\frac{1}{2}\left(\norm{x_1}_{\mathcal{L}_1}^2+\norm{x_2}_{\mathcal{L}_2}^2\right).
	\end{equation}
 Furthermore in this paper we consider
\begin{equation}
\label{defG}
{\mathcal G} =G_0+G_1\frac{\partial}{\partial \zeta},
\end{equation}
	with $G_0, G_1 \in \mathbb{R}^{n\times n}$ and $G_1$ full rank. ${\mathcal G}^*$ is the formal adjoint of ${\mathcal G}$ {\it i.e.}
	\begin{equation*}
		{\mathcal G}^*=G_0^T-G_1^T\frac{\partial}{\partial \zeta}.
	\end{equation*}
	 The dissipation operator $R$ is bounded, symmetric ($R^*=R$) and coercive ($\langle z,Rz\rangle_{L_2}>a \norm{z}_{L_2}$, $\forall z \in L_2(\left[ 0,L\right], \mathbb{R}^n)$ and $a>0$). $u_d$ and $y_d$ denote the distributed input and output, respectively. The system (\ref{twoconservationlaws}-\ref{output}) with \eqref{defG} stems from the modeling
	 of wavelike systems like elastic strings, Timoshenko beams or waves and beams organised in networks. The proposed approach is easy to extend to second order operators defining Euler-Bernouilli beam equation for example. We define $P_1= \begin{bmatrix} 0 & G_1\\ G_1^T &0 \end{bmatrix}.$
	
	\begin{definition}\label{Boundaryvariables}
		The boundary port variables associated to the system \eqref{twoconservationlaws} are defined by: 
		\begin{equation}\label{Bound_port_param}
			\begin{bmatrix}
				f_\partial\\e_\partial
			\end{bmatrix} = \underbrace{\frac{1}{\sqrt{2}}\begin{bmatrix}
				P_1 & -P_1\\I & I
			\end{bmatrix}}_{R_{\text{ext}}}\begin{bmatrix}
				\mathcal{L}x(L)\\
				\mathcal{L}x(0)\\
			\end{bmatrix}.
		\end{equation}	\end{definition}
	
	By definition the boundary port variables are such that:
	\begin{equation} \label{eq:time_deriveH}
	\frac{\dd H}{\dd t} =  \int_0^L y_d^*u_d\dd\zeta+f_\partial^T e_\partial-\int^L_0 (\mathcal{L}_2x_2)^*
	R (\mathcal{L}_2x_2) d\zeta.
	\end{equation}

\begin{thm}\label{theorem:plant} 
	Let $W$ be a $2n\times 4n$ matrix. If $W$ has full rank and satisfies $W\Sigma W^T\geq 0$, where $\Sigma =\begin{bmatrix}
		0& I \\ I & 0
	\end{bmatrix}
	$, then the system operator 
	\[
		\mathcal{A}= \left(\mathcal{J} - \mathcal{R}\right)\mathcal{L}
	\] where \[\mathcal{J}=\begin{bmatrix}0& \mathcal{G}\\-\mathcal{G}^*&0 \end{bmatrix}\;\; \mbox{and}\;\; \mathcal{R}=\begin{bmatrix}0&0 \\ 0 & R\end{bmatrix}\]
	with domain 
	\begin{equation*} 
		D(\mathcal{A} ) = \left\{x \in H^N\left(\left[ 0,L\right], \mathbb{R}^{2n}\right) \mid  \begin{bmatrix}
			f_\partial\\e_\partial
		\end{bmatrix} \in \text{ker} (W) \right\}
	\end{equation*}
	generates a contraction semigroup. 
\end{thm}
\begin{pf}
	The proof follows Theorem 4.1 in \cite{le2005dirac}.
\end{pf}
Boundary inputs and outputs are defined by:
\begin{equation} \label{eq:boundary}
	\begin{aligned}
		u_b&=W\begin{bmatrix}f_{\partial}^T&e_{\partial}^T\end{bmatrix}^T,&y_b&=\tilde{W}\begin{bmatrix}f_{\partial}^T&e_{\partial}^T\end{bmatrix}^T,
	\end{aligned}
\end{equation}
with $\tilde{W}$ full rank and $\begin{bmatrix}
	W^T,&\tilde{W}^T
\end{bmatrix}^T$ invertible. 

With the plant system established, the first step in the design procedure using an early lumping approach is to spatially discretize \eqref{twoconservationlaws}. The discretization needs to preserve the structure and the passivity of the system to take advantage of the PHS properties. Therefore we apply the mixed finite element method \cite{golo2004hamiltonian} and the approximated system of \eqref{twoconservationlaws} is again a PHS with $p$ elements:
\begin{subequations}
	\label{discretizedVS_origin}
	\begin{align}
		\begin{bmatrix}\dot{x}_{1d}\\\dot{x}_{2d}
		\end{bmatrix}&=\left(J_n-R_n\right)\begin{bmatrix}Q_{1}x_{1d}\\Q_{2}x_{2d}\end{bmatrix}+B_bu_b+\begin{bmatrix}
			0\\B_{0d}
		\end{bmatrix}\mathbf{u}_d,\\
		y_b&=B_b^T\begin{bmatrix}
			Q_{1}x_{1d}\\Q_{2}x_{2d}
		\end{bmatrix}+D_bu_b,\\
		\mathbf{y}_d&=\begin{bmatrix}
			0&B_{0d}^T
		\end{bmatrix}\begin{bmatrix}
			Q_{1}x_{1d}\\Q_{2}x_{2d}
		\end{bmatrix},
	\end{align}	
\end{subequations}
where $x_{id}=\begin{bmatrix}x_i^1&\cdots&x_i^p\end{bmatrix}^T$ for $i \in \left\{1,\cdots,2n\right\}$, $\mathbf{u}_d\in\mathbb{R}^p$, $\mathbf{y}_d\in\mathbb{R}^p$,
\begin{align*}
	J_n&=\begin{bmatrix}
		0&J_i\\-J_i^T&0
	\end{bmatrix} &\text{and}&
	& R_n&=\begin{bmatrix}
		0&0\\0&R_d
	\end{bmatrix},
\end{align*} 
are the discretized matrices of the operators $\mathcal{J}$ and $\mathcal{R}$ with $J_i$ and $R_d$ the discretized matrices of the operators $\mathcal{G}$ and $R$. $B_{0d}\in\mathbb{R}^{np\times p}$, $Q_1\in\mathbb{R}^{np\times np}$ and $Q_2\in\mathbb{R}^{np\times np}$ are the discretized matrices of $B_0$, $\mathcal{L}_1$ and $\mathcal{L}_2$, respectively.  
The input $u_b$ denotes the boundary input which corresponds to the boundary actuation or/and conditions. 
Since the distributed actuation of the system is considered,  we assume that there is no energy changes (actuation) at the boundary of the spatial domain, {\it i.e.} $u_b = 0$ and the discretized system \eqref{discretizedVS_origin} can therefore be simplified.

The Hamiltonian of the discretized model \eqref{discretizedVS_origin} writes:
\begin{equation} 
	H_d(x_{1d},x_{2d})=\frac{1}{2}\left(x_{1d}^TQ_{1}x_{1d}+x_{2d}^TQ_{2}x_{2d}\right).
	\label{discretizeHamiltonian}
\end{equation}
It is important to notice that in what follows the choice of the structure preserving discretization method is not unique. One could have alternatively used other discretization methods such as \cite{kotyczka2018weak,moulla2012pseudo} that also guarantee the existence of port Hamiltonian structure and structural invariants suitable for control design purposes. 
Furthermore we consider a finite number of inputs for control design. In this case, the infinite dimensional system \eqref{twoconservationlaws} is in general not controllable but stabilizable, because the uncontrollable modes are already exponentially stable.

	\section{Control by interconnection and energy shaping}
	\label{CIES}
\begin{figure}[htbp!]
	\centering
	\begin{tikzpicture}
		\draw(-.3,0) rectangle (3.5,1);
		\node at (1.6,.5) {\textsf{Discretized plant system}};
		\draw(0,-2.5) rectangle (3.2,-1.5);
		\node at (1.6,-2) {\textsf{Controller}};
		\draw(0.1,-0.75) circle(0.1);
		\draw[->] (0.1,-1.5)--(0.1,-0.85);
		\node at (-.05,-1) {$-$};
		\draw[->] (0.1,-0.65)--(0.1,0);
		\draw[->] (0.3,0)--(0.3,-1.5);
		\node at (0.7,-0.75) {$\cdots$};
		\node at (0.95,-1) {$-$};
		\draw(1.1,-0.75) circle(0.1);
		\draw[->] (1.1,-1.5)--(1.1,-0.85);
		\draw[->] (1.1,-0.65)--(1.1,0);
		\draw[->] (1.3,0)--(1.3,-1.5);
		\node at (1.42,-1) {$-$};
		\draw(1.5,-0.75) circle(0.1);
		\draw[->] (1.5,-1.5)--(1.5,-0.85);
		\draw[->] (1.5,-0.65)--(1.5,0);
		\draw[->] (1.7,0)--(1.7,-1.5);
		\node at (1.82,-1) {$-$};
		\draw(1.9,-0.75) circle(0.1);
		\draw[->] (1.9,-1.5)--(1.9,-0.85);
		\draw[->] (1.9,-0.65)--(1.9,0);
		\draw[->] (2.1,0)--(2.1,-1.5);
		\node at (2.5,-0.75) {$\cdots$};
		\node at (2.75,-1) {$-$};
		\draw(2.9,-0.75) circle(0.1);
		\draw[->] (2.9,-1.5)--(2.9,-0.85);
		\draw[->] (2.9,-0.65)--(2.9,0);
		\draw[->] (3.1,0)--(3.1,-1.5);
		\draw (-2,0.5)--(-0.3,0.5);
		\node at (-1,0.7) {\textsf{BC}$(0)$};
		\draw (3.5,0.5)--(5.2,0.5);
		\node at (4.2,0.7) {\textsf{BC}$(L)$};
		\node at (-.3,-0.4) {\textcolor{blue}{$\textbf{u}_d$}};
		\node at (3.4,-0.4) {\textcolor{blue}{$\textbf{y}_d$}};
		\node at (3.4,-1.1) {\textcolor{blue}{$u_c$}};
		\node at (-0.3,-1.1) {\textcolor{blue}{$y_c$}};
	\end{tikzpicture}
	\caption{Distributed control by interconnection strategy.}
	\label{Fig.ControlScheme}
\end{figure}
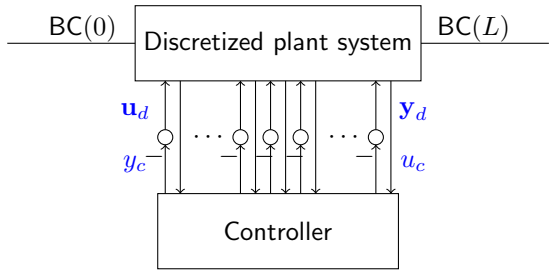

	In this section, we extend the CbI method to the in-domain distributed input and output case. The main difference with CbI for finite dimensional PHSs \cite{van2000l2,ortega2001putting} is that the controller uses the overall information of the plant into consideration, as depicted in Fig.~\ref{Fig.ControlScheme}. As a result, one can shape the distributed Hamiltonian function all over the system with an appropriate parametrization of the controller and the use of structural invariants {\it i.e.} Casimir functions \cite{duindam2009modeling}. The main objective of the proposed CbI method is to improve the closed loop performances over a given frequency range while guaranteeing the overall closed loop stability (v.s. neglected dynamics during the synthesis). One can also modify the equilibrium point by changing the minima of the energy function.

The controller in Fig.~\ref{Fig.ControlScheme} is designed to be a finite dimensional PHS, which is expressed as follows:
\begin{equation} 
	\label{controller}
	\begin{aligned}
		\dot{x}_c&=\left[J_c-R_c\right]Q_cx_c+B_cu_c,\\
		y_c &= B_c^TQ_cx_c+D_cu_c,
	\end{aligned}
\end{equation}
where $x_c \in \mathbb{R}^{m }$, $J_c=-J_c^T \in \mathbb{R}^{m \times m}$, $R_c=R_c^T\geq 0$ and $Q_c=Q_c^T\geq 0$, $B_c\in\mathbb{R}^{m\times m}$, $\mathbb{R}^{m \times m} \ni D_c>0$, $u_c\in\mathbb{R}^{m}$ and $y_c\in\mathbb{R}^{m}$. Matrices $Q_c$ and $D_c$ are used for energy shaping and damping injection/diffusion, respectively. 

Without considering external signals, the interconnection between the discretized plant system \eqref{discretizedVS_origin} and the controller \eqref{controller} is given by
\begin{equation} 
	\label{interconnect}
	\begin{bmatrix}
		\mathbf{u}_d\\
		u_c
	\end{bmatrix}=\begin{bmatrix}
		0&-M\\
		M^T&0
	\end{bmatrix}\begin{bmatrix}
		\mathbf{y}_d\\y_c
	\end{bmatrix},
\end{equation}
where $M=I_m\otimes\mathbf{1}_k\in \mathbb{R}^{p\times m}$, with $\otimes$ denoting the Kronecker product,  $I_m$ representing the identity matrix of dimension $m$ and $k$ is the number of elements covered by one actuator. 

	The passive interconnection \eqref{interconnect} guarantees the passivity of the closed-loop system. It results in a new PHS in closed-loop:
	\begin{equation}
		\label{closeloop}
		\dot{x}_{cl}=\left[J_{cl}-R_{cl}\right]Q_{cl}x_{cl},
	\end{equation}
	where $x_{cl}=\begin{bmatrix}x_{1d}^T,&x_{2d}^T,&x_c^T\end{bmatrix}^T$,  $Q_{cl}=\text{diag}\begin{bmatrix}Q_{1},&Q_{2},&Q_{c}\end{bmatrix}$,
		\begin{align*}
		&J_{cl}=\begin{bmatrix}
			O&J_i&0\\
			-J_i^T&0&-B_{0d}MB_c^T\\
			0&B_cM^TB_{0d}^T&J_c
		\end{bmatrix},\\
		&R_{cl}=\begin{bmatrix}
			0&0&0\\
			0&R_d+B_{0d}MD_cM^TB_{0d}^T&0\\
			0&0&R_c
		\end{bmatrix}.
	\end{align*}
	The Hamiltonian of the controller \eqref{controller} is:
	\begin{equation}
		H_c(x_c)= \frac{1}{2}x_c^TQ_cx_c.
		\label{H_c}
	\end{equation}
	Therefore, the closed-loop Hamiltonian function reads:
	\begin{equation}
		H_{cld}(x_{1d},x_{2d},x_c)=H_d(x_{1d},x_{2d})+H_c(x_c).
		\label{closeloopHamiltonian_origin}
	\end{equation}
 The next step is to design the controller matrices $J_c$, $R_c$, $B_c$, $Q_c$, and $D_c$ in order to shape the closed-loop Hamiltonian $\eqref{closeloopHamiltonian_origin}$. In Proposition \ref{Pro_1} we first show how the controller states can be related to the plant states  using structural invariants.
	\begin{prop1}
	\label{Pro_1}
	Choosing $J_c=0$, and $R_c=0$, the closed-loop system \eqref{closeloop} admits the Casimir function $C(x_{1d},x_c)$ defined by:
	\begin{equation} 
		C(x_{1d},x_c)=B_cM^TB_{0d}^TJ_i^{-1}x_{1d}-x_c
		\label{invariantget}
	\end{equation}
	as structural invariant, {\it i.e.} $\dot{C}(x_{1d},x_c)=0$ along the closed-loop trajectories. If the initial conditions of $x_{1d}(0)$ and $x_c(0)$ satisfy $C(x_{1d}(0),x_c(0))=0$, the controller is a proportional-integral control, and the control law \eqref{controller}-\eqref{interconnect} is equivalent to the state feedback:
	\begin{equation} 
		\begin{aligned}
			y_c&=B_c^TQ_cB_cM^TB_{0d}^TJ_i^{-1}x_{1d}+D_cM^TB_{0d}^TQ_2x_{2d},\\
			\mathbf{u}_d&=-My_c.
		\end{aligned}
		\label{statefeedback}
	\end{equation}
	Therefore, the closed-loop system yields:
	\vspace{-5mm}
	\begin{small}
		\begin{equation} 
			\begin{bmatrix}
				\dot{x}_{1d}\\
				\dot{x}_{2d}
			\end{bmatrix}=\begin{bmatrix}
				0&J_i\\
				-J_i^T&-\tilde{R}_{d}
			\end{bmatrix}\begin{bmatrix}
				\tilde{Q}_{1}x_{1d}\\
				Q_{2}x_{2d}
			\end{bmatrix},
			\label{closenew}
		\end{equation} 
	\end{small}
	where 
	\begin{align} 
	\tilde{R}_{d}=\left[R_d+B_{0d}MD_cM^TB_{0d}^T\right], \\
		\tilde{Q}_{1}=Q_{1}+J_i^{-T}B_{0d}MB_c^TQ_cB_cM^TB_{0d}^TJ_i^{-1}
		\label{eq:tildQ1}
	\end{align}
	are the new closed-loop dissipation matrix and energy matrix associated to $x_{1d}$.
\end{prop1}
	\begin{pf} 
		We consider here Casimir functions of the form:
		\begin{equation}
			C(x_{1d},x_{2d},x_c)=F(x_{1d},x_{2d})-x_c.
		\end{equation}
		The time derivative of $C$ is given by
		\begin{equation} 
		\label{invariant}
		\begin{aligned}
			\frac{\dd C}{\dd t}&=\frac{\partial^T C}{\partial x_{cl}}\frac{\partial x_{cl}}{\partial t}\\
			&=\begin{bmatrix}
				\frac{\partial^T F}{\partial x_{1d}},& \frac{\partial^T F}{\partial x_{2d}}, &-I
			\end{bmatrix}\left(J_{cl}-R_{cl}\right)e_{cl},
		\end{aligned}	
	\end{equation}
	where $e_{cl}=\frac{\partial H_{cld}}{\partial x_{cl}}=Q_{cl} x_{cl}$.
	The Casimir functions are dynamic invariants, i.e. $\dot{C}=0$ that do not depend on the trajectories of the system {\it i.e.} on the Hamiltonian. Therefore, \eqref{invariant} with  $\dot{C}=0$ gives rise to the following matching equations:
	\vspace{-5mm}
	\begin{small}
		\begin{subequations}
			\begin{align}
				&\frac{\partial^T F}{\partial x_{2d}}\left(-J_i^T\right)=0,\label{x2n0}\\
				&\frac{\partial^T F}{\partial x_{1d}}J_i-\frac{\partial^T F}{\partial x_{2d}}\tilde{R}_{d}-B_cM^TB_{0d}^T=0,\label{Fx2}\\
				&\frac{\partial^T F}{\partial x_{2d}}\left(-B_{0d}MB_c^T\right)-\left(J_c-R_c\right)=0 \label{eq:JcRc}.
			\end{align}
	\end{subequations}	\end{small}
	
	Solving \eqref{x2n0}, one gets $\partial F/\partial x_{2d}=0$, which indicates that $x_c$ does not depend on $x_{2d}$. Therefore, with $J_c=-J_c^T$ and $R_c=R_c^T>0$, \eqref{eq:JcRc} indicates that $J_c$ and $R_c$ equal zero. Since $J_i$ is full rank, from \eqref{Fx2} one gets \eqref{invariantget} as a structural invariant as soon as the initial condition $x_c(0)$ is chosen properly. Taking the initial conditions  $x_{1d}(0)$ and $x_c(0)$ such that $C(x_{1d}(0),x_c(0))=0$, \eqref{invariantget} becomes
	\begin{equation} 
		B_cM^TB_{0d}^TJ_i^{-1}x_{1d}-x_c=0,
		\label{eq:change}
	\end{equation}
	which allows to link the state of the controller with the state of the plant. Replacing $x_c$ in \eqref{closeloop} by \eqref{eq:change}, the control law \eqref{controller} becomes a state feedback as formulated in \eqref{statefeedback}. Therefore the closed-loop system \eqref{closeloop} becomes \eqref{closenew}.
\end{pf}
	From Proposition \ref{Pro_1}, the closed-loop Hamiltonian function \eqref{closeloopHamiltonian_origin} is now only function of the discretized plant state variables :
	\begin{equation}
		H_{cld}(x_{1d},x_{2d})=\frac{1}{2}\left(x_{1d}^T\tilde{Q}_{1}x_{1d}+x_{2d}^TQ_{2}x_{2d}\right),
		\label{closeloopHamiltonian}
	\end{equation}
	with its time derivative being:
\begin{equation} 
	\frac{\dd H_{cld}}{\dd t}=-x_{2d}^TQ_{2}\left(R_d+B_{0d}MD_cM^TB_{0d}^T\right)Q_{2}x_{2d} \leq 0.
	\label{timederivHam}
\end{equation}

From a physical point of view, \eqref{closeloopHamiltonian} implies that with the dynamic controller \eqref{controller} equivalent to the state feedback \eqref{statefeedback}, it is possible to change, at least partially (depending on $p$ and the range of $B_{0d}$), the energy matrix related to $x_{1d}$. Actually, one can only shape the energy matrix related to the first $p$ elements of $x_{1d}$, i.e. $\left(\tilde{Q}_1\right)_{p\times p}$. For a given number of distributed input $m$, the objectives of the energy shaping is to look for matrices $B_c$ and $Q_c$ such that the norm of the difference (considered here in the Frobenius norm, see Definition 6.4 of \cite{shores2007applied}) between the desired energy matrix $\tilde{Q}_{1d}$ and the closed loop one $\tilde{Q}_1$ is minimal: 
\begin{equation} 
	\label{eq:pro1_orig}
	\underset{
		B_c^TQ_cB_c
	}{\text{min}}
	\norm{J_i^{-T}B_{0d}MB_c^TQ_cB_cM^TB_{0d}^TJ_i^{-1}+Q_{1}-\tilde{Q}_{1d} }_F.
\end{equation}  
If we consider $p$ elements and eliminate $B_{0d}$, \eqref{eq:pro1_orig} is equivalent to:
\begin{equation} 
	\label{eq:pro1_refor}
	\underset{
		B_c^TQ_cB_c
	}{\text{min}}
	\norm{\left(J_i\right)_{p\times p}^{-T}MB_c^TQ_cB_cM^T\left(J_i\right)_{p\times p}^{-1}-Q_m}_F,
\end{equation}
where the $\left(J_i\right)_{p\times p}$ are the first $p$ lines $p$ columns of $J_i$ and $Q_m=\left(\tilde{Q}_{1d}-Q_{1}\right)_{p\times p}\geq 0$. Furthermore, \eqref{eq:pro1_refor} can be formalized by the optimization Problem \ref{Problem}.
	\begin{problem}\label{Problem}
		The closed loop energy function related to the first $p$ elements of $x_{1d}$  is shaped in an optimal way if and only if  $X=B_c^TQ_cB_c\in SR_0^{m\times m}$  minimizes the criterion 
		\begin{equation}\label{ls}
			f(X)=\norm{AXA^T-Q_m }_F,
		\end{equation}
		where $A=\left(J_i\right)_{p\times p}^{-T}M\in \mathbb{R}^{p\times m}$ and
	$SR_0^{m\times m}$ represents the set of symmetric and positive semi-definite matrices. 
\end{problem}
	The solution to Problem \ref{Problem} depends on the independent number of distributed input that are available. We consider two different cases: the {\it ideal} fully-actuated case ($m=p$) and the under-actuated case ($m<p$).
	\subsection{Fully-actuated case}\label{secfullactuated}
\begin{figure}[htbp!]
		\centering
		\includegraphics[scale=0.5]{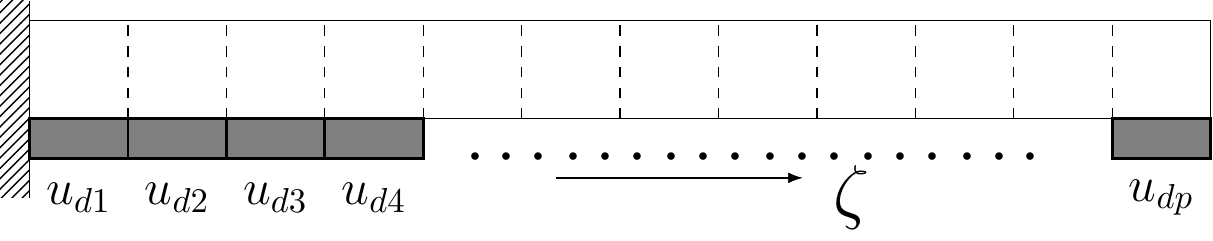}
		\caption{Fully-actuated case illustration.}
		\label{fullyactuated}
	\end{figure}

	We first consider the fully-actuated case where each discretized element of the plant is controlled by an independent input, {\it i.e.} $u_d\in \mathbb{R}^m$ and $m=p$, as illustrated in Fig.~\ref{fullyactuated}. 
The input matrix $M=I\in\mathbb{R}^{p\times p}$. Therefore, the optimization Problem \ref{Problem} admits an exact solution that is given in Proposition \ref{pro:fully}.
\begin{prop1}
	\label{pro:fully}
	In the fully-actuated case, {\it i.e.} $m=p$ the optimization Problem \ref{Problem} has an exact analytical solution $\hat{X}=M^{-1}\left(J_i\right)_{p\times p}^T Q_m \left(J_i\right)_{p\times p}M^{-T}$ leading to $f(X)=0$. The controller matrices $B_c$ and $Q_c$ can be chosen as:
	\begin{equation} 
		\label{fullydesign}
		\begin{aligned}
			B_c&=\left(J_i\right)_{p\times p},&Q_c&= Q_m.
		\end{aligned}
	\end{equation}
\end{prop1}
\begin{pf}
	The matrix $A$ is invertible, therefore, \eqref{ls} admits a minimum in $0$ when 
	\begin{equation} 
		\label{xmesure}
		\begin{aligned}
			\hat{X}&=A^{-1}Q_mA^{-T}=M^{-1}\left(J_i\right)_{p\times p}^TQ_m\left(J_i\right)_{p\times p}M^{-T}\\
			&=\left(J_i\right)_{p\times p}^TQ_m\left(J_i\right)_{p\times p}.
		\end{aligned}
	\end{equation}
	From the expression of $X$, one can choose $B_c$ and $Q_c$ as in \eqref{fullydesign} to satisfy \eqref{xmesure}.
\end{pf}
\begin{remark}
	The choice $B_c=\left(J_i\right)_{p\times p}$ can be regarded as the numerical approximation of the operator $\frac{\partial}{\partial\zeta}$, which has also been used in the late lumping control design approach in \cite{tranchantcdc}.
\end{remark}
	
	\subsection{Under-actuated case}
	\label{secOne2elements}
\begin{figure}[thpb!] 
		\vspace{6pt}
		\centering
		\includegraphics[width=5cm]{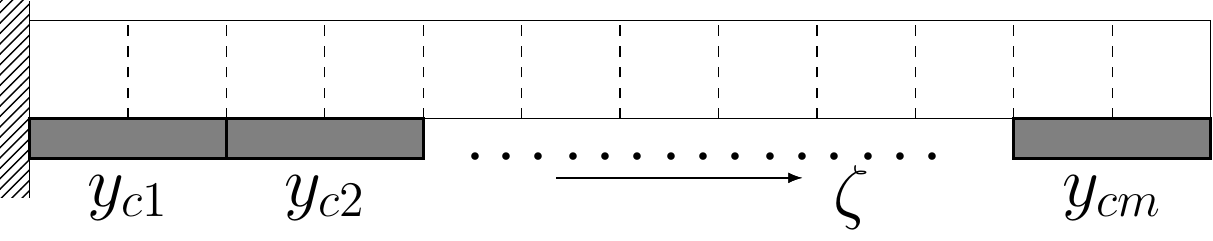}
		\caption{Under-actuated case with $k=2$.}
		\label{sameInput}
	\end{figure}

	We study now the more realistic case where the same control input is applied to a set of elements, as shown in Fig.~\ref{sameInput}, where $k$ denotes the number of elements sharing the same input. The number of distributed inputs $m$ is less than the number of discretized elements, and follows $m=p/k$.

	In this case the controller has less degree of freedom than in the fully actuated case, hence the matrix $A$ in \eqref{ls} is not invertible and the optimization Problem \ref{Problem} is ill-conditioned.
The solution of the optimization Problem \ref{Problem} is given in Proposition \ref{Pro_3}.

\begin{prop1}\label{Pro_3}
	$f(X)$ defined in \eqref{ls} is convex and the minimization of $f(X)$ is equivalent to the minimization of $f^2(X)$, which has a unique minimum given for $\hat{X}=V\Sigma_0^{-1}U_1^TQ_mU_1\Sigma_0^{-1}V^T$, with $V$, $\Sigma_0$ and $U_1$ the matrices of the singular value decomposition (SVD) of the matrix A {\it i.e.}
	\begin{align}
		A=U\Sigma V^T=\begin{bmatrix}
			U_1&U_2
		\end{bmatrix}\begin{bmatrix}
			\Sigma_0\\0
		\end{bmatrix}V^T,
		\label{svd}
	\end{align}
	where $U\in\mathbb{R}^{p\times p}$ and $V\in\mathbb{R}^{m\times m}$ are unitary matrices, $U_1\in\mathbb{R}^{p\times m}$, $U_2\in\mathbb{R}^{p\times q}$, $q=p-m$, and $\Sigma_0=\Sigma_0^T\geq 0$ is the diagonal matrix of singular values of $A$.
\end{prop1}
\begin{pf}	
	The proof of Proposition \ref{Pro_3} is similar to that of Proposition 3 in \cite{Liu2021}.
	Substituting \eqref{svd} into $f^2(X)$, one gets:
	\begin{equation} 
		\begin{aligned}
			&\underset{X\in SR_0^{m\times m}}{\text{min}}f^2(X)\\
			&=\underset{X\in SR_0^{m\times m}}{\text{min}}\norm{U \Sigma V^TXV\Sigma^TU^T-Q_m }_F^2\\
			&=\underset{X\in SR_0^{m\times m}}{\text{min}}\left(\norm{\Sigma_0V^TXV\Sigma_0^T-T_1}_F^2+2\norm{T_2}_F^2+\norm{T_3}_F^2\right),
		\end{aligned}
		\label{svdfsquare}
	\end{equation}
	where $T_1=U_1^TQ_mU_1$, $T_2=U_1^TQ_mU_2$, and $T_3=U_2^TQ_mU_2$. Since $\norm{T_2}_F^2$ and $\norm{T_3}_F^2$ are given once the matrices $A$ and $Q_m$ are defined, the minimization of \eqref{svdfsquare} is equivalent to:
	\begin{equation} 
		\underset{\bar{X}\in SR_0^{m\times m}}{\text{min}}\norm{\bar{X}-T_1}_F^2,\ \text{with} \ \bar{X}=\Sigma_0V^TXV\Sigma_0^T.
		\label{PT}
	\end{equation}
	According to Theorem 2.1 in \cite{highamLAA1988}, $T_1 \in SR_0^{m\times m}$, and \eqref{PT} admits a unique solution $\hat{\bar{X}}=T_1$. Therefore, \eqref{svdfsquare} has the minimum when:
	\begin{equation} 
		\hat{X}=V\Sigma_0^{-1}\hat{\bar{X}}\Sigma_0^{-1}V^T=V\Sigma_0^{-1}U_1^TQ_mU_1\Sigma_0^{-1}V^T.
		\label{xoptimal_under}
	\end{equation}
\end{pf}
The choice of controller matrices $B_c$ and $Q_c$ is not unique, as long as they satisfy the condition \eqref{xoptimal_under}. We will present a possible choice in Subsection \ref{subsection:simu_under}.
\begin{remark}
	\label{rem:damping_injection}
	We have investigated the choices of controller matrices $B_c$ and $Q_c$ under two different cases in Proposition \ref{pro:fully} and \ref{Pro_3}, respectively. The objective is to shape the closed loop Hamiltonian $H_{cld}$ with the modification of part of the potential energy matrix. The choice of the controller matrix $D_c$ follows the similar procedure, with the optimization of the difference between \eqref{timederivHam} and the desired one. 
\end{remark}
	\subsection{Closed-loop stability}
\begin{figure}[htbp!]
	\centering
	\begin{tikzpicture}
		\draw(0,0) rectangle (3.2,1);
		\node at (1.6,.5) {\textsf{Plant system}};
		\draw(0,-2.5) rectangle (3.2,-1.5);
		\node at (1.6,-2) {\textsf{Controller}};
		\draw(0.1,-0.75) rectangle (0.65, -0.85);
		\draw[->] (0.3,-1.5)--(0.3,-0.85);
		\node at (.15,-1) {$-$};
		\draw[->] (0.1,-0.85)--(0.1,0);
		\draw[->] (0.4,-0.85)--(0.4,-1.5);
		\draw[->] (0.2,0)--(0.2,-0.75);
		\draw[->] (0.3,-0.75)--(0.3,0);
		\draw[->] (0.4,0)--(0.4,-0.75);
		\draw[->] (0.5,-0.75)--(0.5,0);
		\draw[->] (0.6,0)--(0.6,-0.75);
		\draw[->] (0.95,-0.85)--(0.95,-1.5);
		\draw(0.65,-0.75) rectangle (1.2, -0.85);
		\draw[->] (0.85,-1.5)--(0.85,-0.85);
		\node at (.7,-1) {$-$};
		\draw[->] (0.65,-0.85)--(0.65,0);
		\draw[->] (0.75,0)--(0.75,-0.75);
		\draw[->] (0.85,-0.75)--(0.85,0);
		\draw[->] (0.95,0)--(0.95,-0.75);
		\draw[->] (1.05,-0.75)--(1.05,0);
		\draw[->] (1.15,0)--(1.15,-0.75);
		\node at (1.65,-0.75) {$\cdots$}; 
		\draw(2.1,-0.75) rectangle (2.65, -0.85);
		\draw[->] (2.3,-1.5)--(2.3,-0.85);
		\draw[->] (2.4,-0.85)--(2.4,-1.5);
		\node at (2.15,-1) {$-$};
		\draw[->] (2.1,-0.85)--(2.1,0);
		\draw[->] (2.2,0)--(2.2,-0.75);
		\draw[->] (2.3,-0.75)--(2.3,0);
		\draw[->] (2.4,0)--(2.4,-0.75);
		\draw[->] (2.5,-0.75)--(2.5,0);
		\draw[->] (2.6,0)--(2.6,-0.75);
		\draw(2.65,-0.75) rectangle (3.2, -0.85);
		\draw[->] (2.85,-1.5)--(2.85,-0.85);
		\node at (2.7,-1) {$-$};
		\draw[->] (2.65,-0.85)--(2.65,0);
		\draw[->] (2.75,0)--(2.75,-0.75);
		\draw[->] (2.85,-0.75)--(2.85,0);
		\draw[->] (2.95,0)--(2.95,-0.75);
		\draw[->] (3.05,-0.75)--(3.05,0);
		\draw[->] (3.15,0)--(3.15,-0.75);
		\draw[->] (2.95,-0.85)--(2.95,-1.5);
		\draw (-2,0.5)--(-0,0.5);
		\node at (-1,0.7) {\textsf{BC}$(0)$};
		\draw (3.2,0.5)--(5.2,0.5);
		\node at (4.2,0.7) {\textsf{BC}$(L)$};
		\node at (-.3,-0.4) {\textcolor{blue}{$u_d$}};
		\node at (3.4,-0.4) {\textcolor{blue}{$y_d$}};
		\node at (3.4,-1.1) {\textcolor{blue}{$u_c$}};
		\node at (-0.3,-1.1) {\textcolor{blue}{$y_c$}};
	\end{tikzpicture}
	\caption{Distributed control by interconnection strategy.}
	\label{Fig.ControlSchemeUnder}
\end{figure}
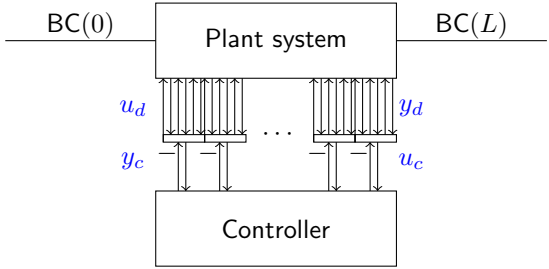

	In this subsection we consider the closed-loop stability of the infinite-dimensional system \eqref{twoconservationlaws} controlled by the finite-dimensional controller \eqref{controller} derived from the early lumping approach.  The power-preserving interconnection between \eqref{twoconservationlaws} and \eqref{controller} is formulated as:
\begin{equation} 
	\begin{bmatrix}
		u_d\\u_c
	\end{bmatrix}=\begin{bmatrix}
		0&-\mathbf{1}_\zeta\\
		\mathbf{1}_\zeta^*&0
	\end{bmatrix}\begin{bmatrix}
		y_d\\y_c
	\end{bmatrix},
	\label{eq:interconnection}
\end{equation}
with
\begin{align}
	\mathbf{1}_\zeta&:\mathbb{R}^m\rightarrow L_2,&\mathbf{1}_\zeta^*&:L_2\rightarrow\mathbb{R}^m.
\end{align}
$\mathbf{1}_\zeta$ is the characteristic function that distributes the point-wise value of the controller in $R^{m}$ space to the sub-interval $L_2$ space, as illustrated in Fig. \ref{Fig.ControlSchemeUnder}.

\begin{lemma}
	The interconnection \eqref{eq:interconnection} generates a Dirac structure with the following power conservation:
	\begin{equation}
		\int_{0}^{L}y_d^*\mathbf{1}_\zeta y_c\dd\zeta=y_c^T\mathbf{1}_\zeta^*y_d.
		\label{eq:powerconservation}
	\end{equation}
\end{lemma}
The closed-loop system is equivalent to:
\begin{equation} 
	\dot{\mathcal{X}}=\underbrace{\begin{bmatrix}
			0&\mathcal{G}&0\\
			-\mathcal{G}^*&-\mathcal{R}_{cl}&-B_0\mathbf{1}_\zeta B_c^T\\
			0&B_c\mathbf{1}_\zeta^*B_0^*&0
		\end{bmatrix}\mathcal{L}_{cl}
	}_{\mathcal{A}_{cl}}\mathcal{X},
	\label{eq:closed}
\end{equation}
where $\mathcal{X}=\begin{bmatrix}
	x\\x_c
\end{bmatrix}\in X_s$ is the state defined on the state space $X_s= L_2\left([0,L],\mathbb{R}^{2n} \right)\times \mathbb{R}^m$, $\mathcal{R}_{cl}=R+B_0\mathbf{1}_\zeta D_c\mathbf{1}_\zeta^*B_0^*$ and $\mathcal{L}_{cl}=\text{diag}\left(\mathcal{L}_1,\mathcal{L}_2,Q_c\right)$.

The Hamiltonian of \eqref{eq:closed} is:
\begin{equation} 
	H_{cl}=\frac{1}{2}\left(\norm{x_1}_{\mathcal{L}_1}^2+\norm{x_2}_{\mathcal{L}_2}^2\right)+\frac{1}{2}x_c^TQ_cx_c
	\label{eq:Hamiltonianclosed}
\end{equation}
with \vspace{-5mm}\begin{small} \begin{equation} 
		\begin{aligned}
			\frac{\dd H_{cl}}{\dd t}&=\int_0^L  y_d^*u_d\dd\zeta+y_c^Tu_c-\int^L_0 (\mathcal{L}_2x_2)^*\mathcal{R}_{cl} (\mathcal{L}_2x_2) d\zeta\\
			&=-\int^L_0 (\mathcal{L}_2x_2)^*\mathcal{R}_{cl} (\mathcal{L}_2x_2) d\zeta.
		\end{aligned}
		\label{eq:TimeHcl}
	\end{equation}
\end{small}
The last step of \eqref{eq:TimeHcl} is derived considering \eqref{eq:powerconservation}.

	In order to prove stability of the closed loop system using Lyapunov arguments and LaSalle's invariance principle we first propose the following theorems. 
	\begin{thm} \label{theorem_compact}
		The linear operator $\mathcal{A}_{cl}$ defined in \eqref{eq:closed} generates a contraction semigroup on $X_s$.
	\end{thm}
	\begin{pf}
		To prove that the closed loop operator $\mathcal{A}_{cl}$ generates a contraction semigroup,  we apply Lumer-Phillips Theorem (Theorem 1.2.3 in \cite{liu1999semigroups}). The proof is done in two steps: first, we show that the operator $\mathcal{A}_{cl}$ is dissipative. Second, we show that \begin{equation}
		\label{eq:range}
		\text{range}\left(\lambda I-\mathcal{A}_{cl}\right)\in X_s, \text{ for }\lambda>0.
		\end{equation}  
		
		According to Definition 6.1.4 in \cite{Jacob2012}, $\mathcal{A}_{cl}$ is dissipative if Re$\langle \mathcal{A}_{cl}\mathcal{X}, \mathcal{X}\rangle\leq 0$, which is equivalent to $\langle \mathcal{A}_{cl}\mathcal{X}, \mathcal{X}\rangle+\langle \mathcal{X},\mathcal{A}_{cl}\mathcal{X}\rangle\leq 0$. For the sake of clarity and without any restriction, we take $\mathcal{L}_1=\mathcal{L}_2=1$ and $Q_c=I$ in the rest of this section. From \eqref{eq:closed}, one has:
		\begin{equation}
		\label{eq:dissipativeA}
			\begin{aligned}
				&\langle \mathcal{A}_{cl}\mathcal{X}, \mathcal{X}\rangle+\langle \mathcal{X},\mathcal{A}_{cl}\mathcal{X}\rangle\\
			=&\langle \mathcal{G}x_2,x_1\rangle_{L_2}+\langle -\mathcal{G}^*x_1-\mathcal{R}_{cl}x_2,x_2\rangle_{L_2}\\
			&+\langle -B_0\mathbf{1}_\zeta B_c^Tx_c,x_2\rangle_{L_2}+\langle B_c\mathbf{1}_\zeta^*B_0^*x_2,x_c\rangle_{ \mathbb{R}^m}\\
			&+\langle x_1, \mathcal{G}x_2\rangle_{L_2}+\langle x_2, -\mathcal{G}^*x_1-\mathcal{R}_{cl}x_2\rangle_{L_2}\\
			&+\langle x_2,-B_0\mathbf{1}_\zeta B_c^Tx_c\rangle_{L_2}+\langle x_c,B_c\mathbf{1}_\zeta^*B_0^*x_2\rangle_{ \mathbb{R}^m}.
			\end{aligned}
		\end{equation}
	According to \eqref{eq:powerconservation}, we get 
	\begin{equation}
		\begin{aligned}
			\langle B_0\mathbf{1}_\zeta B_c^Tx_c,x_2\rangle_{L_2}&=\langle B_c\mathbf{1}_\zeta^*B_0^*x_2,x_c\rangle_{ \mathbb{R}^m},\\
			\langle x_2,-B_0\mathbf{1}_\zeta B_c^Tx_c\rangle_{L_2}&=\langle x_c,B_c\mathbf{1}_\zeta^*B_0^*x_2\rangle_{ \mathbb{R}^m}.
		\end{aligned}
	\label{eq:interconnectInnerProduct}
	\end{equation}
		Substituting \eqref{eq:interconnectInnerProduct} into \eqref{eq:dissipativeA}, we have:
		\begin{equation*}
			\begin{aligned}
				&\langle \mathcal{A}_{cl}\mathcal{X}, \mathcal{X}\rangle+\langle \mathcal{X},\mathcal{A}_{cl}\mathcal{X}\rangle\\
				=&\langle \mathcal{G}x_2,x_1\rangle_{L_2}+\langle -\mathcal{G}^*x_1,x_2\rangle_{L_2}-\langle\mathcal{R}_{cl}x_2,x_2\rangle_{L_2}\\
				&+\langle x_1, \mathcal{G}x_2\rangle_{L_2}+\langle x_2, -\mathcal{G}^*x_1\rangle_{L_2}-\langle x_2,-\mathcal{R}_{cl}x_2\rangle_{L_2}\\
				=&-\langle\mathcal{R}_{cl}x_2,x_2\rangle_{L_2}-\langle x_2,-\mathcal{R}_{cl}x_2\rangle_{L_2}\leq 0,
			\end{aligned}
		\end{equation*}
		where the last step is obtained according to the boundary conditions. Therefore, the operator $\mathcal{A}_{cl}$ is dissipative.
		
		To show \eqref{eq:range}, we apply the proof of Theorem 3.3.6 in \cite{Augner2016} with adjustment dedicated to our in-domain control. For the sake of simplicity, we choose $\lambda=1$. Taking an arbitrary function 
		\begin{equation*}
			f=\begin{bmatrix}
				\tilde{x}\\\tilde{x}_c
			\end{bmatrix}\in X_s,
		\end{equation*}
	 \eqref{eq:range} is then equivalent to the problem:
		\begin{equation}
			\text{find } \mathcal{X}\in X_s: \left(I-\mathcal{A}_{cl}\right)\mathcal{X}=f,
			\label{eq:rangeProblem}
		\end{equation}
		which is again equivalent to:
		\begin{subequations}
			\begin{align}
				x_1-\mathcal{G}x_2&=\tilde{x}_1,  \label{eq:ode1}\\
				\mathcal{G}^*x_1+\left(I+\mathcal{R}_{cl}\right)x_2+B_0\mathbf{1}_\zeta B_c^Tx_c&=\tilde{x}_2, \label{eq:ode2}\\
				-B_c\mathbf{1}_\zeta^*B_0^*x_2+x_c&=\tilde{x}_c. \label{eq:ode3}
			\end{align}
		\end{subequations}
		Substituting \eqref{eq:ode3} into \eqref{eq:ode2}, one gets:
		\begin{equation}
			\mathcal{G}^*x_1+\left(I+\mathcal{M}\right)x_2=\tilde{x}_2-B_0\mathbf{1}_\zeta B_c^T\tilde{x}_c,
			\label{eq:ode5}
		\end{equation}
	with $\mathcal{M}=\mathcal{R}_{cl}+B_0\mathbf{1}_\zeta B_c^T\mathbf{1}_\zeta^*B_0^*$.
	
	According to the definition of $\mathcal{G}$ and $\mathcal{G}^*$, \eqref{eq:ode1} and \eqref{eq:ode5} become:
		\begin{align}
				x_1-\sum_{k=0}^{1}G_k\frac{\partial^k}{\partial\zeta^k}x_2&=\tilde{x}_1,  \label{eq:odeG1}\\
				\sum_{k=0}^{1}(-1)^{k}G_k^T\frac{\partial^k}{\partial\zeta^k}x_1+\left(I+\mathcal{M}\right)x_2&=\tilde{x}_f,\label{eq:odeG2}
		\end{align}
	with $\tilde{x}_f=\tilde{x}_2-B_0\mathbf{1}_\zeta B_c^T\tilde{x}_c$.  Thus one gets:
	\begin{small}
	\begin{align}
		\frac{\partial^N x_1}{\partial \zeta^N}&=(-1)^{N+1}G_N^{-T}\left[\sum_{k=0}^{N-1}(-1)^kG_k^T\frac{\partial^k x_1}{\partial \zeta^k}+\left(I+\mathcal{M}\right)x_2-\tilde{x}_f\right],\label{eq:odeGNx1}\\
			\frac{\partial^N x_2}{\partial \zeta^N}&=G_N^{-1}\left[x_1-\sum_{k=0}^{N-1}G_k\frac{\partial^k x_2}{\partial \zeta^k}-\tilde{x}_1\right].\label{eq:odeGNx2}
	\end{align}
	\end{small}
	Define:
		\begin{equation}
			h=\begin{bmatrix}
				x_1\\x_2
			\end{bmatrix}.
		\end{equation}
		According to \eqref{eq:odeGNx1} and \eqref{eq:odeGNx2}, one can derive the following relation:
			\begin{equation}
			\frac{\partial h}{\partial\zeta}=B_hh+g_h.
			\label{eq:relationh}
		\end{equation}
	For $N=1$, 
	\begin{align*}
		B_h&=\begin{bmatrix}
			G_1^{-T}\begin{bmatrix}
				G_0^T,&I+\mathcal{M}
			\end{bmatrix}\\
		G_1^{-1}\begin{bmatrix}
			1,&-G_0
		\end{bmatrix}
		\end{bmatrix},
	&g_h&=\begin{bmatrix}
		-G_1^{-T}\tilde{x}_f\\
		-G_1^{-1}\tilde{x}_1
	\end{bmatrix}.
	\end{align*}
			The solution of the function \eqref{eq:relationh} is derived as:
		\begin{equation}
			h(\zeta)=e^{B_h\zeta}h(0)+q(\zeta),
			\label{eq:solutionh}
		\end{equation}
		with $q(\zeta)=\int_{0}^{\zeta}e^{\zeta-s}B_hg_h\dd s$.
		
		Therefore, to solve the problem \eqref{eq:rangeProblem}, one needs to find the solution of $h(\zeta)$, and eventually the solution of $h(0)$. According to the boundary condition in Theorem 1, we have:
		\begin{equation*}
			\begin{aligned}
				WR_{\text{ext}}\begin{bmatrix}
					h(L)\\h(0)
				\end{bmatrix}&=WR_{\text{ext}}\begin{bmatrix}
					e^{B_h}h(0)+q(L)\\h(0)
				\end{bmatrix}
				=\begin{bmatrix}
					0\\0
				\end{bmatrix}.
			\end{aligned}		
		\end{equation*}
		By calculation, $WR_{\text{ext}}\begin{bmatrix}
			e^{B_h}\\I
		\end{bmatrix}$ has full rank. Hence, one can get the solution of $h(0)$ as:
		\begin{equation*}
			h(0)=-\left(WR_{\text{ext}}\begin{bmatrix}
				e^{B_h}\\I
			\end{bmatrix}\right)^{-1}WR_{\text{ext}}\begin{bmatrix}
			q(L)\\0
		\end{bmatrix}.
		\end{equation*}
	Therefore $h(\zeta)$ is obtained from \eqref{eq:solutionh}. One then has the solution of $x=\begin{bmatrix}
			1&0&\cdots&0
		\end{bmatrix}h(\zeta)$. Substituting $x$ into \eqref{eq:ode2}, one obtains $x_c$. As a result, the problem \eqref{eq:rangeProblem} is solved. According to the Lumer-Phillips theorem, the operator $\mathcal{A}_{cl}$ generates a contraction semigroup that concludes the proof.
	\end{pf}
	\begin{thm}\label{theorem:compact}
		The operator $\mathcal{A}_{cl}$ has a compact resolvent.
	\end{thm}
	\begin{pf}
		According to the Definition A.4.24 in \cite{Curtain1995}, we need to prove that the operator $\left(\lambda I-\mathcal{A}_{cl}\right)^{-1}$ is compact for some $\lambda\in\rho\left(\mathcal{A}_{cl}\right)$, with $\rho\left(\mathcal{A}_{cl}\right)$ denoting the resolvent set of $\mathcal{A}_{cl}$. This proof follows from Garding's inequality (Theorem 7.6.4 in \cite{naylor1971linear}) and the proof of Theorem 2.26 in \cite{Villegas2007}.
		
		Define $\mathcal{T}=\lambda I-\mathcal{A}_{cl}$. From the previous Theorem \ref{theorem_compact}, $\mathcal{A}_{cl}$ generates a contraction semigroup, thus $\lambda>0$ is in the resolvent set of $\mathcal{A}_{cl}$. $\mathcal{T}$ is boundedly invertible and satisfies $\lVert \mathcal{T}\mathcal{X}\rVert_{L_2}\geq \lVert\mathcal{X}\rVert_{H^N}$. Therefore, $\mathcal{T}^{-1}$ is compact which concludes the proof.
%
%
	\end{pf}
	Due to Theorem \ref{theorem_compact} and Theorem \ref{theorem:compact}, the trajectory of the closed-loop system is pre-compact and its asymptotic stability can be proven by Lyapunov arguments and LaSalle's invariance principle (Theorem 3.64 of \cite{luo2012stability}) as shown in Theorem \ref{theorem:stability}.
\begin{thm}
		\label{theorem:stability}
		For any $\mathcal{X}(0)\in L_2\left([0,L],\mathbb{R}^{2n} \right)\times\mathbb{R}^m$, the unique solution of \eqref{eq:closed} tends to zero asymptotically, and the closed-loop system \eqref{eq:closed} is globally asymptotically stable.
	\end{thm}
	\begin{pf}
		We choose the energy of the closed-loop system as Lyapunov function. From \eqref{eq:TimeHcl}, the time derivation of the Lyapunov function is semi-negative definite:
		\begin{equation}
				\frac{\dd H_{cl}}{\dd t}=-\int^L_0 (\mathcal{L}_2x_2)^*\mathcal{R}_{cl} (\mathcal{L}_2x_2) d\zeta\leq 0.
		\end{equation}
		Using LaSalle’s invariance principle, it remains to show that the only solutions associated with $\frac{\dd H_{cl}}{\dd t}$ is $0$ i.e the only solutions associated with $\mathcal{L}_2x_2=0$ is $x_2=0$. Due to the internal dissipation and zero boundary input, the only solution associated with this problem is $0$. The controller being a simple integrator, if well initialized it also converges to $x_c=0$ as the state of the system converges to $x=0$.
	\end{pf}
	
	\section{Numerical simulations}
	\label{SecSimulation}
	As illustrative example we consider a vibrating string of length $L=2\;\si{m}$, modulus of elasticity $T =1.4\times 10^6\;\si{\newton}$ density $\rho=1.225 \; \si{kg/m}$ and dissipation coefficient $R=10^{-3}$. The dynamic model of the string can be written:
	\begin{equation} \label{PHS_vibrating}
	\begin{bmatrix}
		\dot{x}_1\\\dot{x}_2
	\end{bmatrix} = \begin{bmatrix}
		0 & \frac{\partial}{\partial\zeta}\\ \frac{\partial}{\partial\zeta} & -R
	\end{bmatrix}\begin{bmatrix}
		\mathcal{L}_1x_1\\\mathcal{L}_2x_2
	\end{bmatrix} + \begin{bmatrix}
		0\\1
	\end{bmatrix} u_d
\end{equation}
with $x_1(\zeta,t)=\frac{\partial\omega}{\partial\zeta}(\zeta,t)$, and $x_2(\zeta,t)=\rho(\zeta)\frac{\partial\omega}{\partial t}(\zeta,t)$. $\omega(\zeta,t)$ is the longitudinal displacement over the spatial domain with the state space $x=\in L_2(\left[ 0,L \right], \mathbb{R}^2)$.  $\mathcal{L}_1=T$ and $\mathcal{L}_2=\frac{1}{\rho}$. The dissipation term is chosen to be very small $R=10^{-3}$. The distributed input $u_d$ is the force density. From Definition \ref{Boundaryvariables}:
\begin{align*}
	Q&=\begin{bmatrix}
		P_1 & P_2\\-P_2 &0
	\end{bmatrix},
	\text{ with}&P_1&=\begin{bmatrix}
		0 & 1\\1 & 0
	\end{bmatrix},\text{ and}&P_2&=\begin{bmatrix}
		0 & 0\\0 & 0
	\end{bmatrix}. 
\end{align*}
Taking the full rank matrix 
\begin{equation*}
	M_0 = \begin{bmatrix}
		1& 0& 0&0\\0 & 1 &0 &0
	\end{bmatrix}^T,
\end{equation*}
one can compute $Q_M = P_1$ and  $M_Q =  M_0^T$. Then the boundary port variables give:
\begin{equation*}
	\begin{bmatrix}
		f_\partial \\
		e_\partial
	\end{bmatrix} 
	=
	\frac{1}{\sqrt{2}}
	\begin{bmatrix}
		\mathcal{L}_2x_2(L)-\mathcal{L}_2x_2(0)\\
		-\mathcal{L}_1x_1(L)+\mathcal{L}_1x_1(0)\\
		\mathcal{L}_1x_1(L)+\mathcal{L}_1x_1(0)\\
		\mathcal{L}_2x_2(L)+\mathcal{L}_2x_2(0)
	\end{bmatrix}.
\end{equation*}
We consider a clamped-free scenario with in-domain control. Hence, we define the boundary input formulated in \eqref{eq:boundary} with: 
\begin{equation*}
	\begin{aligned}
		W&= \frac{\sqrt{2}}{2} \begin{bmatrix}
			0 & 1 & 1 & 0\\
			-1 & 0 & 0 & 1
		\end{bmatrix},
		&\text{ and } W\Sigma W^T &\geq 0.
	\end{aligned}
\end{equation*}
The clamped-free boundary condition implies $u_b = 0$. 

The discretization matrices in \eqref{discretizedVS_origin} are:
\begin{align*}
	J_i&=\begin{bmatrix}
		\frac{1}{\gamma}\\
		-\frac{1}{\gamma^2}&\frac{1}{\gamma}\\
		\vdots&\ddots&\ddots\\
		(-1)^{p-1}\frac{\left(\gamma'\right)^{p-2}}{\gamma^p}&\cdots&-\frac{1}{\gamma^2}&\frac{1}{\gamma}
	\end{bmatrix}_{p\times p},
\end{align*}
$Q_{1}=\text{diag}\left(T_{ab}\right)\in\mathbb{R}^{p\times p}$, $Q_{2}=\text{diag}\left(\displaystyle \frac{1}{\rho_{ab}}\right)\in\mathbb{R}^{p\times p}$, $R_d=\text{diag}\left(R_{ab}\right)\in\mathbb{R}^{p\times p}$,
with $T_{ab}$, $\displaystyle \frac{1}{\rho_{ab}}$ and $R_{ab}$ chosen to be $\displaystyle \frac{T}{L_{ab}}$, $\displaystyle \frac{1}{\rho L_{ab}}$ and $RL_{ab}$, respectively. $L_{ab}=L/p$. $\gamma$ denotes the effort mapping parameter \cite{golo2004hamiltonian} and $\gamma'=1-\gamma$. They are chosen to be $\displaystyle \frac{1}{2}$ in order to get a centered scheme.

\begin{figure}[thpb]
	\centering
	\includegraphics[scale=0.4]{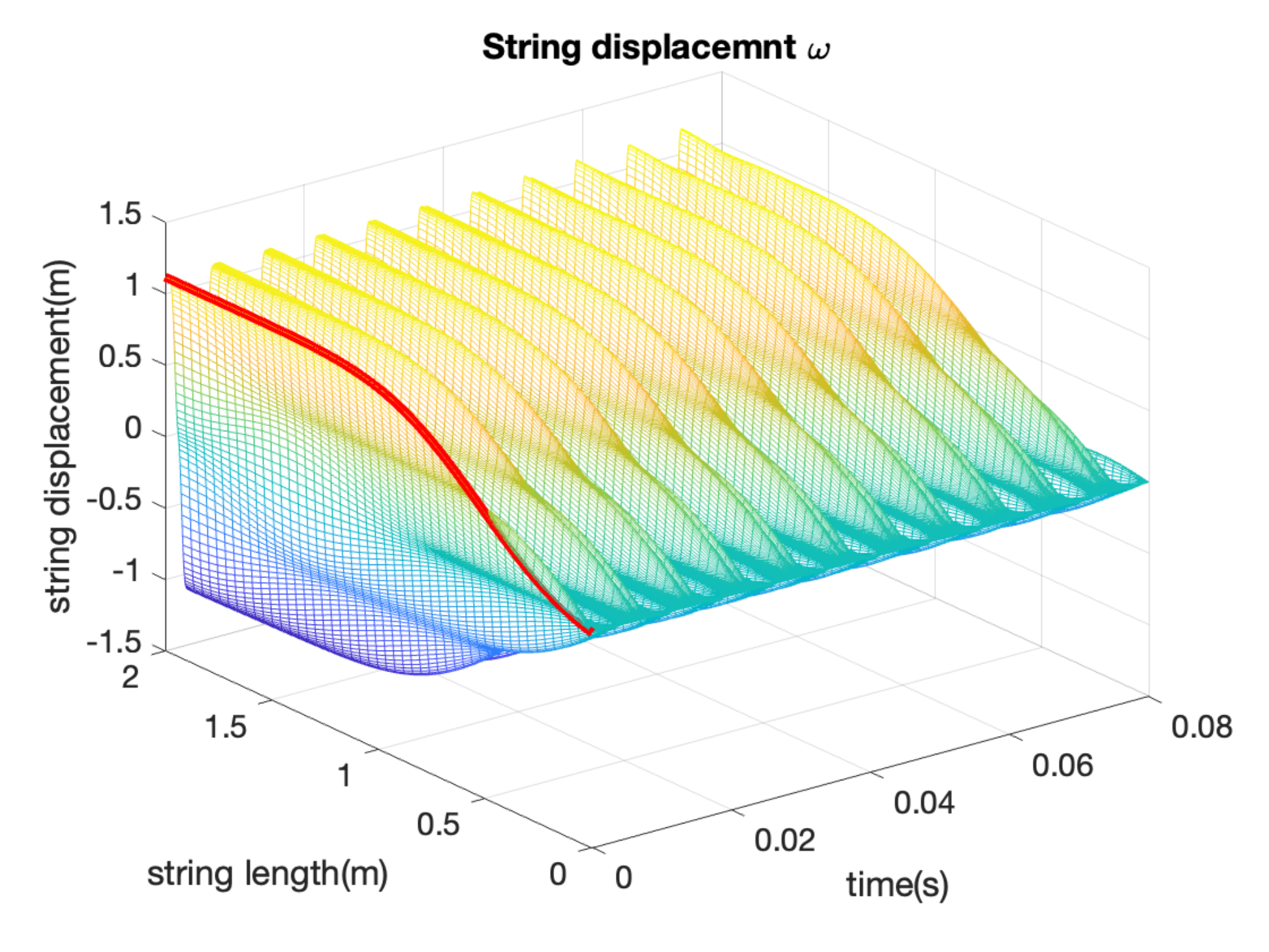}
	\caption{Open loop deformation of the vibrating string.}
	\label{opensimu}
\end{figure}

Initial conditions are set to a spatial distribution $x_1(\zeta,0)\sim\mathcal{N}(1.5,0.113)$ for the strain distribution and to zero for the velocity distribution  {\it i.e.,} $x_2(\zeta,0)=0$. The string is discretized into $50$ elements.
We consider a time step of $5\times 10^{-5}\si{s}$ and mid-point time discretization method\footnote{Implicit midpoint rule is known to be a structure preserving time integrator for PHSs \cite{aoues2013canonical}. It is a particular case in the family of symplectic collocation methods for time integration which is investigated in \cite{kotyczka2018weak}.} for simulations. The open loop evolution of the string deformation $\omega$ is given in Fig.~\ref{opensimu}.

Next we investigate the numerical simulations of the closed-loop system considering both fully-actuated and under-actuated cases.

\subsection{Fully-actuated case}
Following Proposition \ref{Pro_1} and Proposition \ref{pro:fully}, we choose $B_c=J_i$ in order to guarantee the existence of structural invariants, and the initial conditions of the controller such that $C=0$. In this case \eqref{invariantget} becomes:
$x_c=x_{1d}$, and the closed-loop system \eqref{closenew} reads:
\begin{equation} 
	\begin{bmatrix}
		\dot{x}_{1d}\\
		\dot{x}_{2d}
	\end{bmatrix}=\begin{bmatrix}
		0&J_i\\
		-J_i^T&-\left(R_d+D_c\right)
	\end{bmatrix}\begin{bmatrix}
		\tilde{Q}_{1}x_{1d}\\
		Q_{2}x_{2d}
	\end{bmatrix},
\end{equation}
One can see that \emph{the equivalent} closed-loop stiffness ${\tilde{Q}_{1}}$ can be shaped through the choice of  $Q_c$. 

\begin{figure}[thpb]
	\centering
	\subfloat[]{
		\includegraphics[scale=0.28]{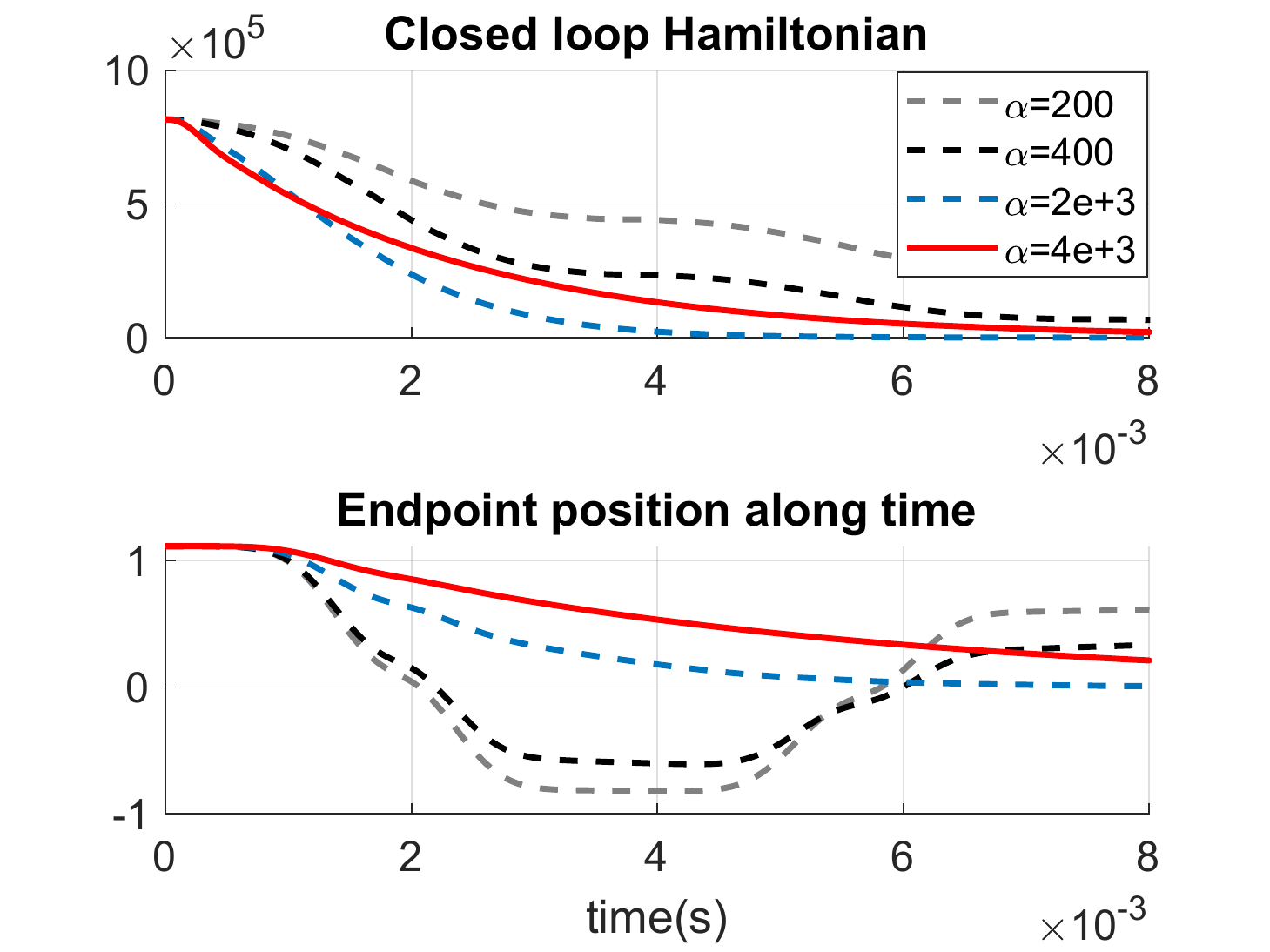}
	}
	\subfloat[]{
		\includegraphics[scale=0.28]{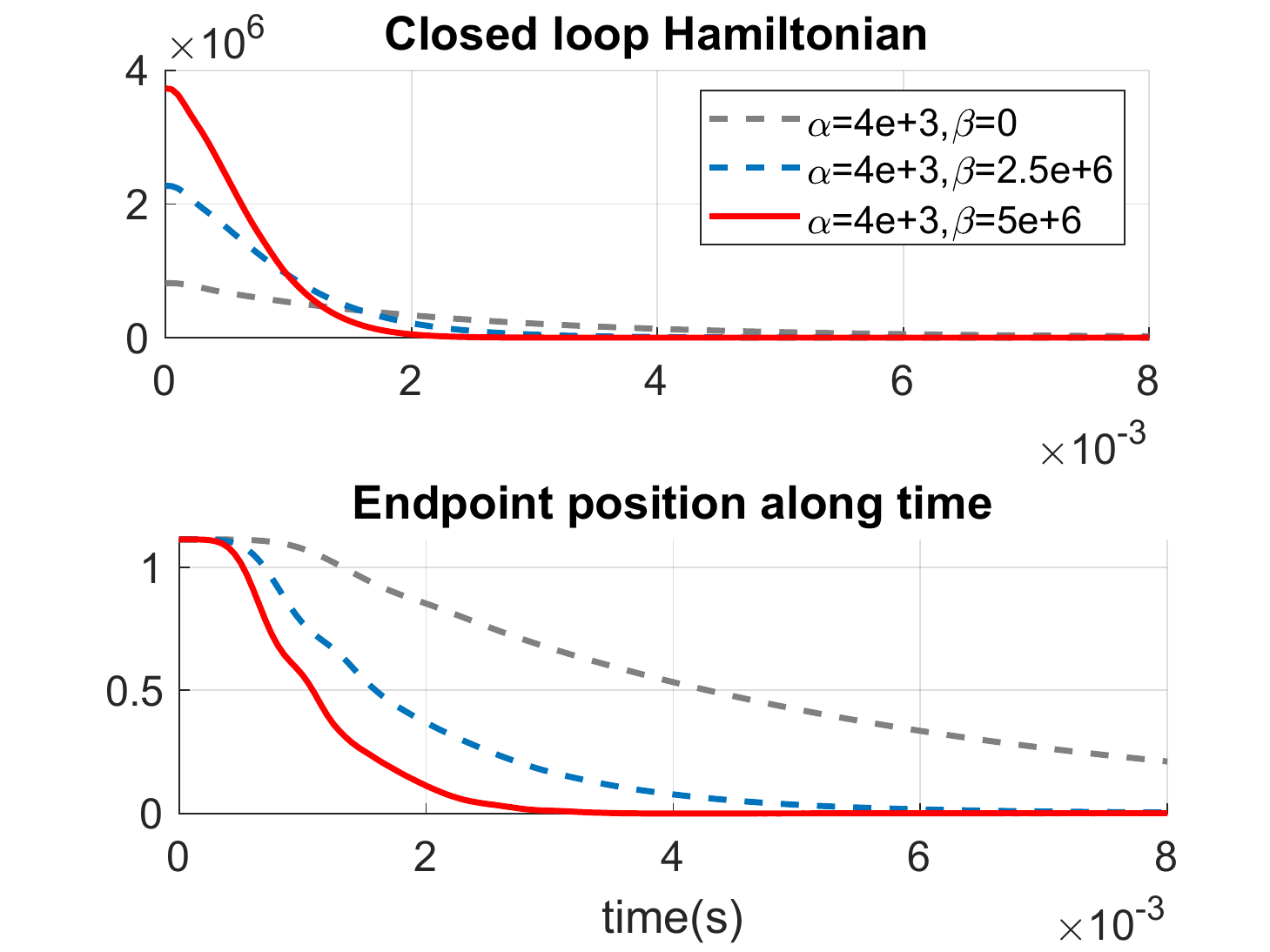}
	}
	\caption{
		Closed-loop Hamiltonian function and endpoint position in the fully-actuated case with (a) pure damping injection and with (b) energy shaping plus damping injection.
	}
	\label{fullyCBIdamping_position}
\end{figure}

We first consider the pure damping injection case, {\it i.e.} varying $D_c$ with $Q_c=0$. We consider $D_c=\text{diag}\left(\alpha L_{ab}\right)$ with $\alpha$ denoting the damping coefficient. In Fig.~\ref{fullyCBIdamping_position}(a) we can see that this degree of freedom allows to damp the vibrations of the string to the detriment of the time response. 

Next we fix $\alpha=4000$ corresponding to the slightly over-damped case in order to illustrate the effect of the energy shaping on the achievable performances. For that we use the control by energy shaping.We can see in Fig.~\ref{fullyCBIdamping_position}(b) that we can speed up the closed-loop system by increasing the closed-loop stiffness via energy shaping, without introducing any overshoot. The energy matrix of the controller $Q_c=\text{diag}\left(\frac{\beta}{L_{ab}}\right)$, with $\beta$ denoting the energy shaping parameter. A good dynamic performance is achieved when $\beta=5\times 10^6$, which relates to an equivalent string stiffness of $ \tilde{T}=6.4\times 10^6$\si{\newton}.

\begin{figure}[thpb]
	\centering
	\subfloat[]{\includegraphics[width=4.3cm]{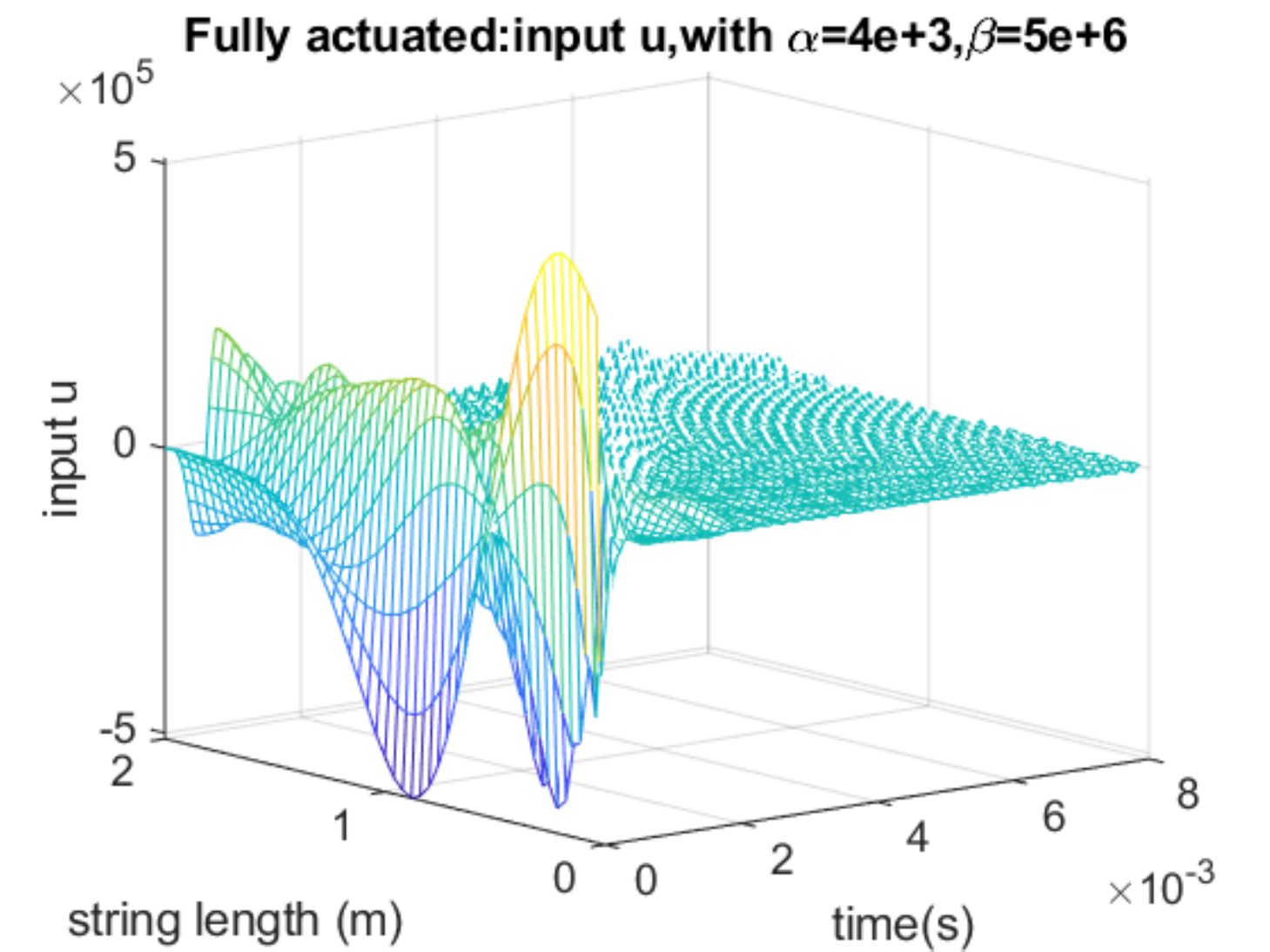}
	}
	\subfloat[]{\includegraphics[width=4.3cm]{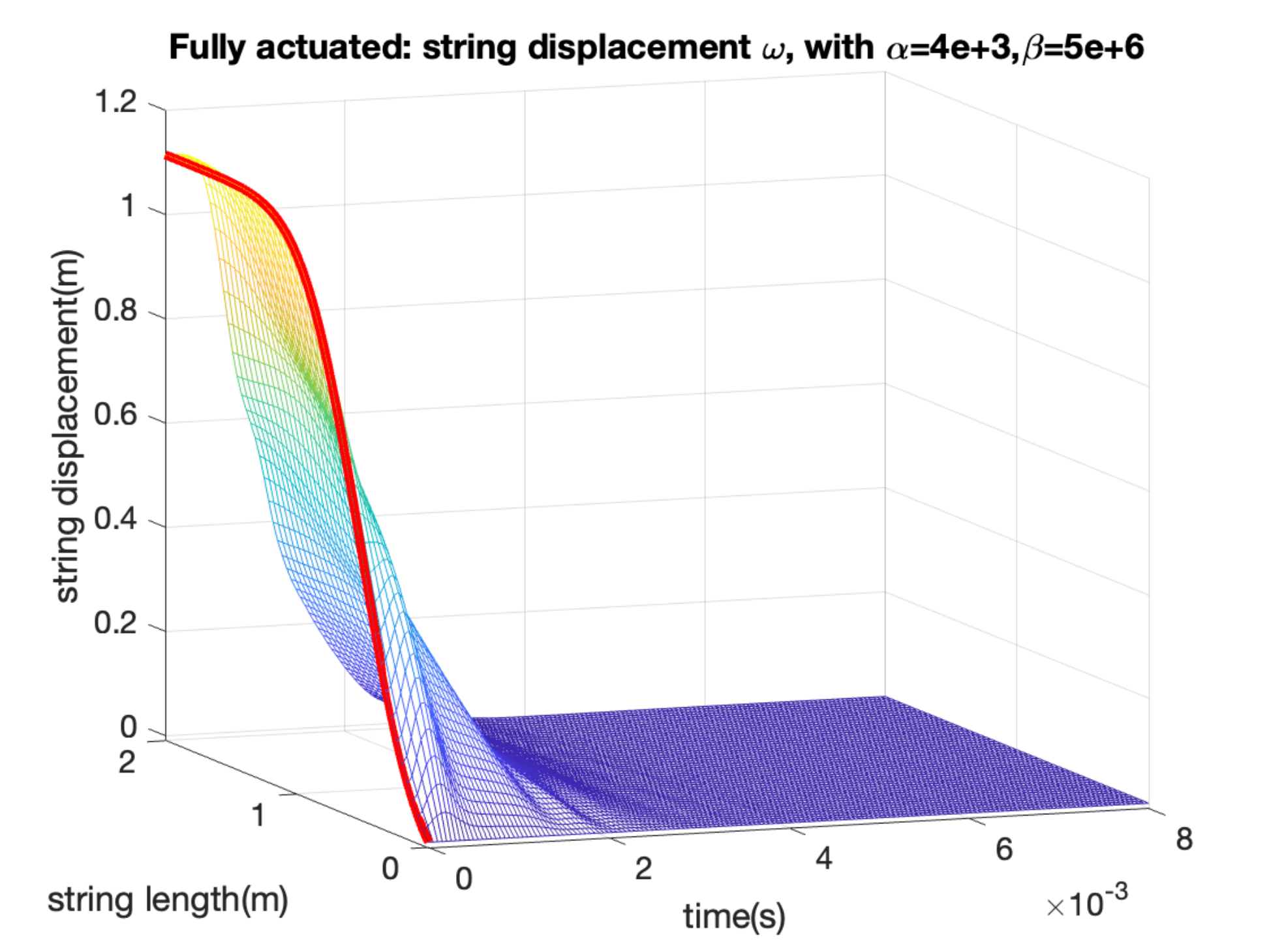}
	}
	\caption{(a) Evolution of the closed-loop input signal and (b) deformation in the energy shaping and damping injection case with full actuation, 
		$\alpha=4\times 10^3$, $\beta=5\times 10^6$.
	}
	\label{fullyCBIdampingstrain}
\end{figure}

The evolution of the distributed input and of the string deformation along time with damping injection and energy shaping are given in Fig.~\ref{fullyCBIdampingstrain}(a) and (b) respectively. We can see in Fig.~\ref{fullyCBIdampingstrain}(a) that the control remains smooth. Fig.~\ref{fullyCBIdampingstrain}(b) shows that the closed-loop stabilization time is about $3\times 10^{-3} s$ which is much faster than $8\times 10^{-3} s$ resulting from the pure damping injection case. 
\subsection{Under-actuated case}
\label{subsection:simu_under}
We now consider that the control is achieved using $m$ patches as depicted in Fig.~\ref{sameInput}. The aim of the control design is to modify as far as possible the internal elasticity $\tilde{T}$ of the string to get similar performances as in the fully-actuated case. We choose the controller matrix $B_c=J_m$ with $J_m\in\mathbb{R}^{m\times m}$ stemming from the discretization of $\frac{\partial}{\partial\zeta}$. According to \eqref{xoptimal_under} in Proposition \ref{Pro_3}, $Q_c=J_m^{-T}V\Sigma_0^{-1}U_1^TQ_mU_1\Sigma_0^{-1}V^TJ_m^{-1}$. 

$D_c$ is chosen according to Remark \ref{rem:damping_injection}, with desired time derivative of the Hamiltonian formulated in \eqref{timederivHam} being the fully-actuated case, {\it i.e.}  in order to satisfy $\underset{D_c\in \mathbb{R}^{m\times m}}{\text{min}}\norm{MD_cM^T- \text{diag}\left(\alpha L_{ab}\right)}_F$. As a results, the optimal $D_c$ is given by $\hat{D}_c=\text{diag}\left(\frac{\alpha L_{ab}}{k}\right)$.

\begin{figure}[thpb]
	\centering
	\subfloat[]{
		\includegraphics[width=4.3cm]{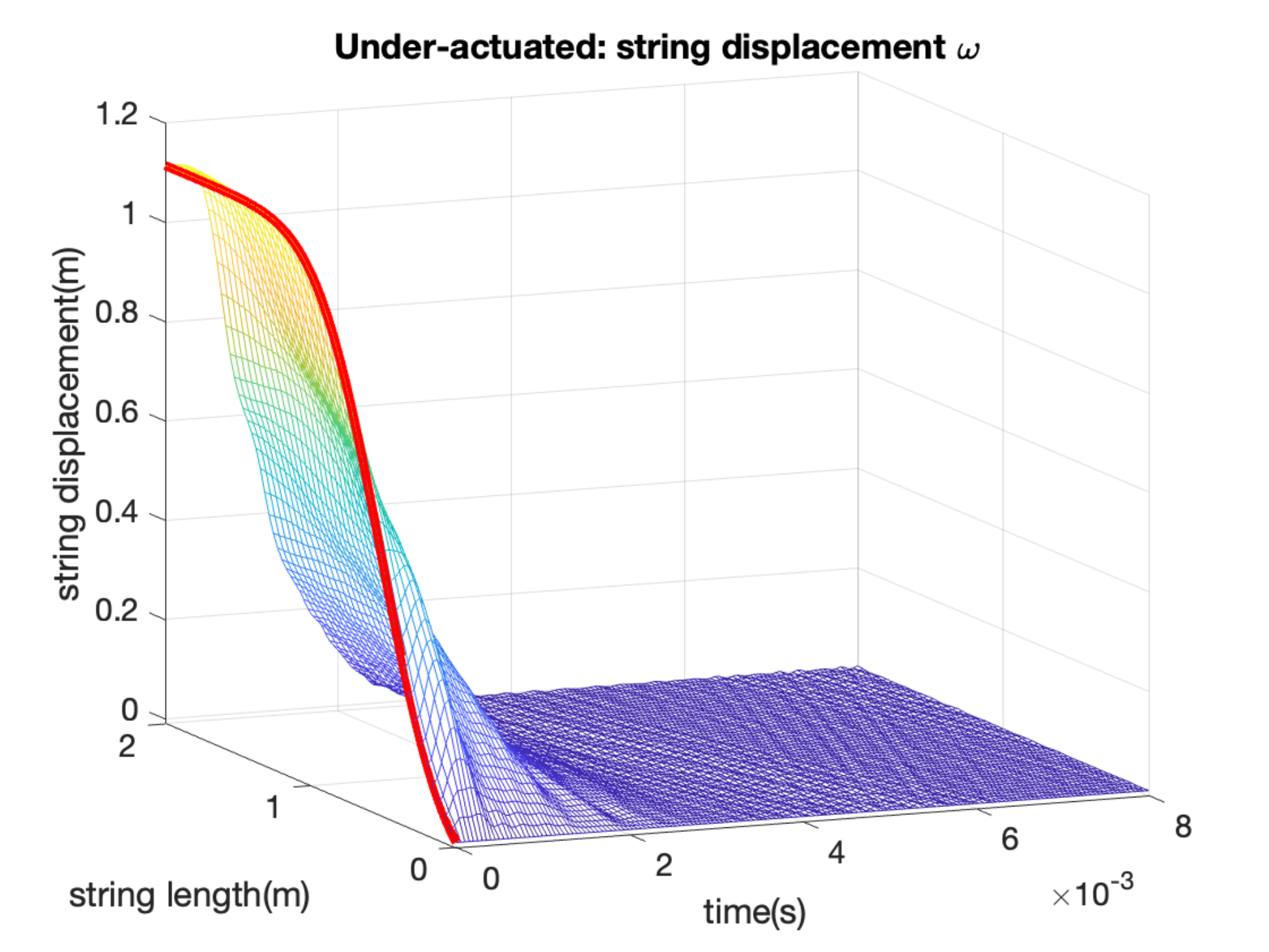}
	}
	\subfloat[]{\includegraphics[width=4.3cm]{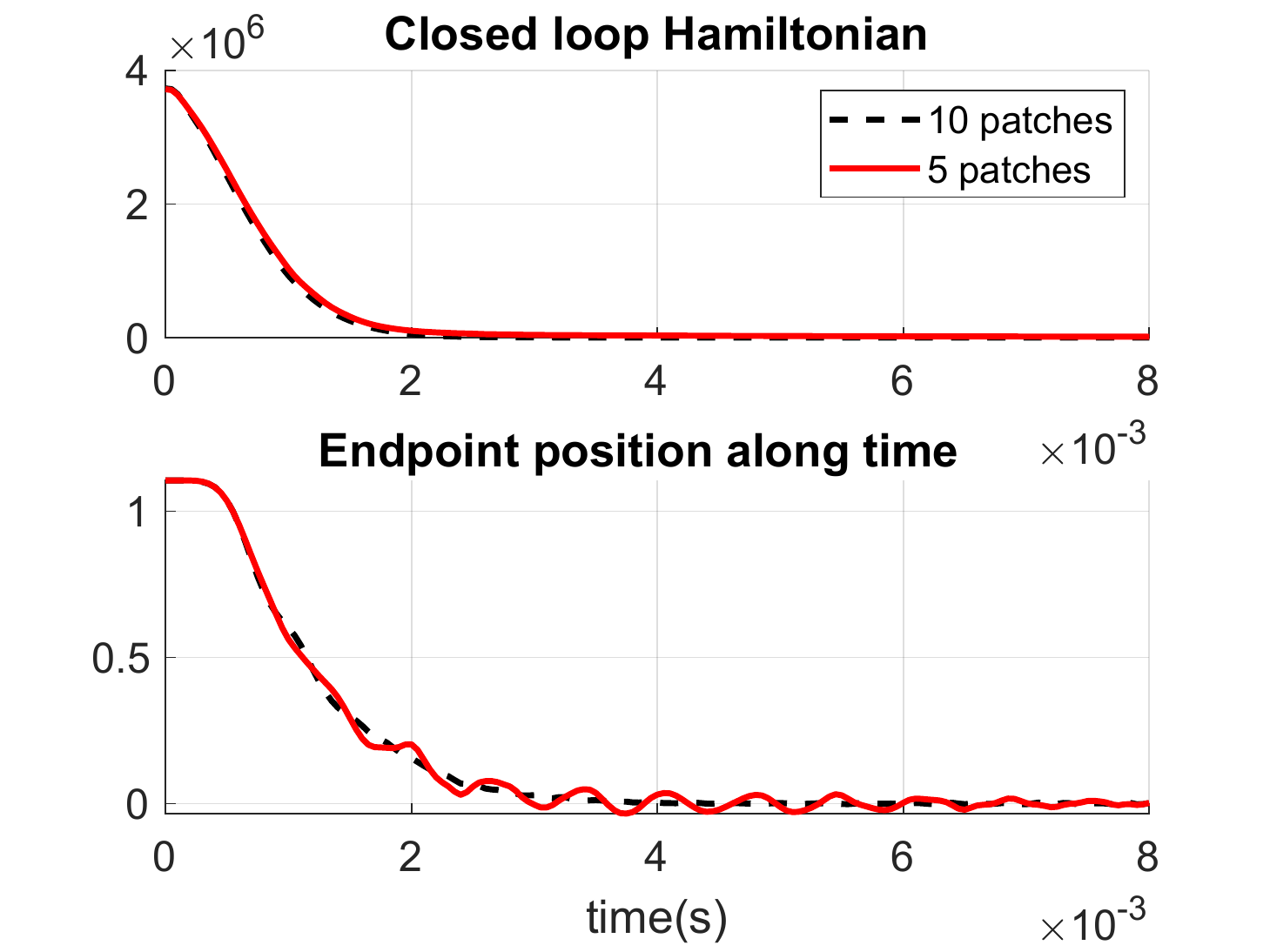}
	}
	\caption{(a) Closed-loop evolution of the deformation, (b) Hamiltonian function and endpoint position in the under-actuated case for $k=5$, and $k=10$.} 
	\label{underactuateSim}
\end{figure}

We first consider the case with 10 patches, {\it i.e.} $p=50$, $m=10$ and $k=5$. The evolution of the string deformation as depicted in Fig.~\ref{underactuateSim}(a) is quite similar to that obtained in the fully-actuated case in Fig.~\ref{fullyCBIdampingstrain}(b). This indicates that if the controller matrices $B_c$, $Q_c$ and $D_c$ are adequately selected, the achievable performances in the under-actuated case can be optimized in order to be close to the ones obtained in the fully-actuated case. When the number of patches is reduced to $5$, these performances are slightly deteriorated at high frequencies as shown in Fig.~\ref{underactuateSim}(b).

\begin{figure}[thpb]
	\centering
	\subfloat[]{\includegraphics[width=4.3cm]{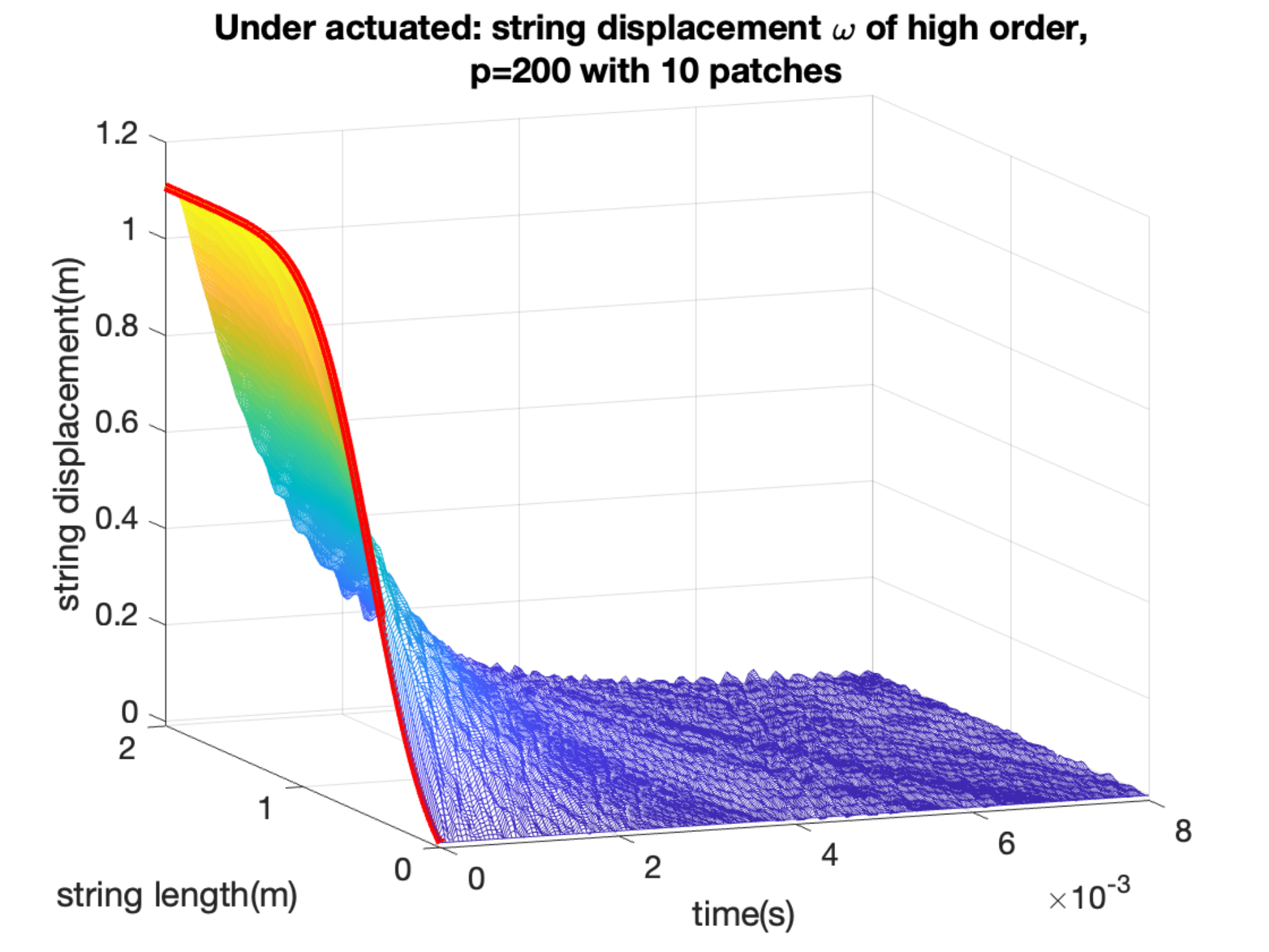}
	}
	\subfloat[]{
		\includegraphics[width=4.3cm]{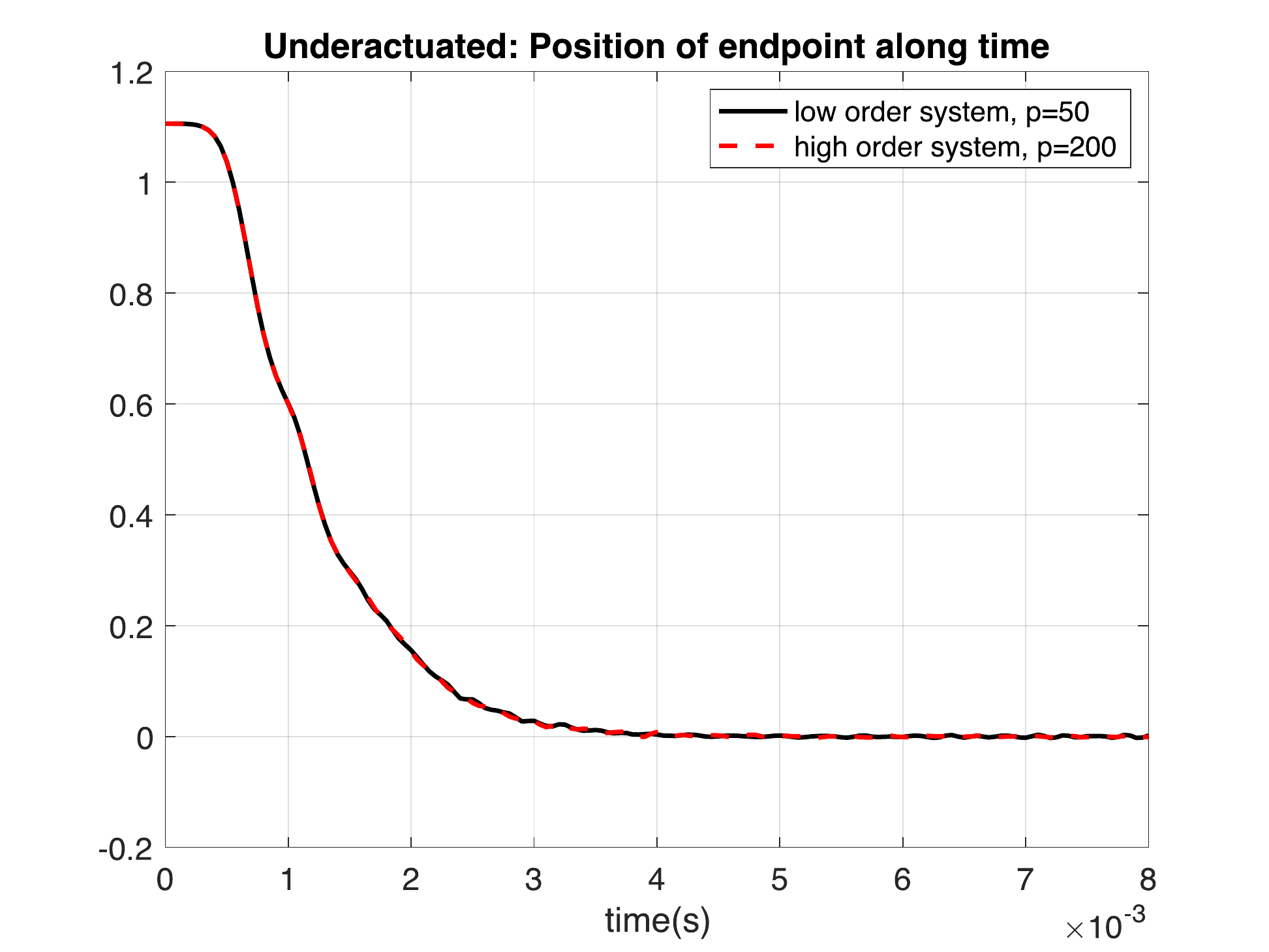}
	}
	\caption{Closed-loop evolution of the deformation of the high order system, and comparison of the endpoint position of the low order and high order systems using the same controller.}
	\label{ValidPrecise}
\end{figure}

In order to illustrate the effect of the neglected dynamics on the achievable performances, we implement the controller designed considering $10$ patches on the discretized system with $p=50$, to a more precise model of the string derived using $p=200$. In Fig.~\ref{ValidPrecise} we can see that, due to the damping injection and the associated closed-loop bandwidth, the neglected dynamics does not impact significantly the closed-loop response of the system to the considered initial condition. An example of in-domain controller design on Timoshenko beam model is investigated in \cite{Liu2021,LeGorrec2022}.

	\section{Frequency analysis of closed-loop systems}	
\begin{figure}[htbp!]
			\centering
			\includegraphics[scale=.4]{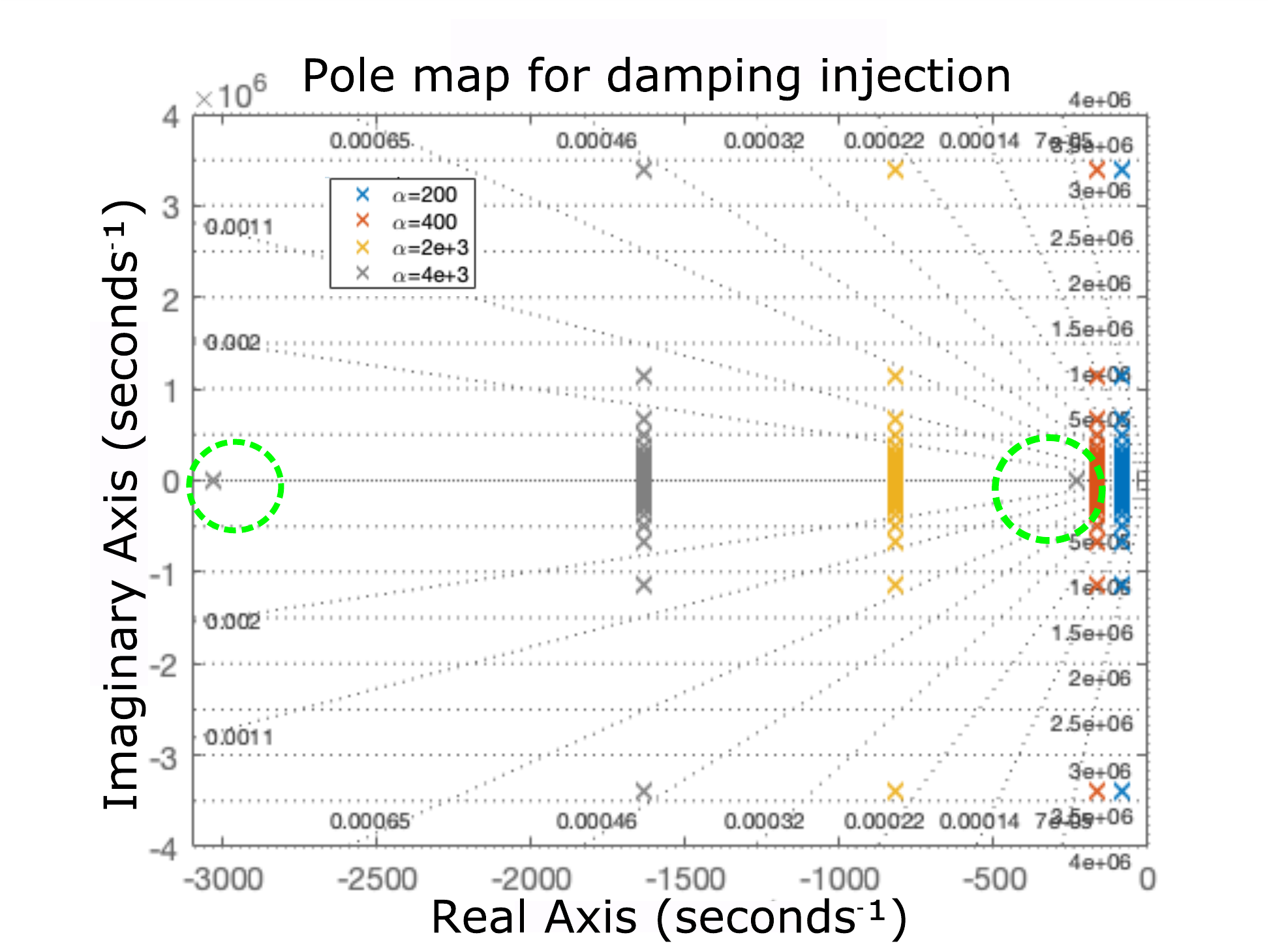}
			\caption{Poles of the closed-loop system with damping injection.}
			\label{fig:PolesZerosDI}
		\end{figure}

In this section, we focus on the frequency analysis of the system. 
		We start with the poles of the closed-loop system with only damping injection actuated on every discretized element of the string corresponding to the simulation results in Fig. \ref{fullyCBIdamping_position}(a). The poles and zero map is illustrated in Fig. \ref{fig:PolesZerosDI}. 
	The damping injection puts the poles away from the imaginary axis in the left hand side of the plot, and increases the stability margin of the residual modes. Moreover, one notice that for $\alpha=4000$, two poles of the closed-loop system are distributed on both sides of the vertical axis. One of the poles located at $-232$ explains the over-damped phenomenon, which is discussed in Fig. \ref{fullyCBIdamping_position}(a).
	
\begin{figure}[htbp!]
			\centering
			\includegraphics[scale=.4]{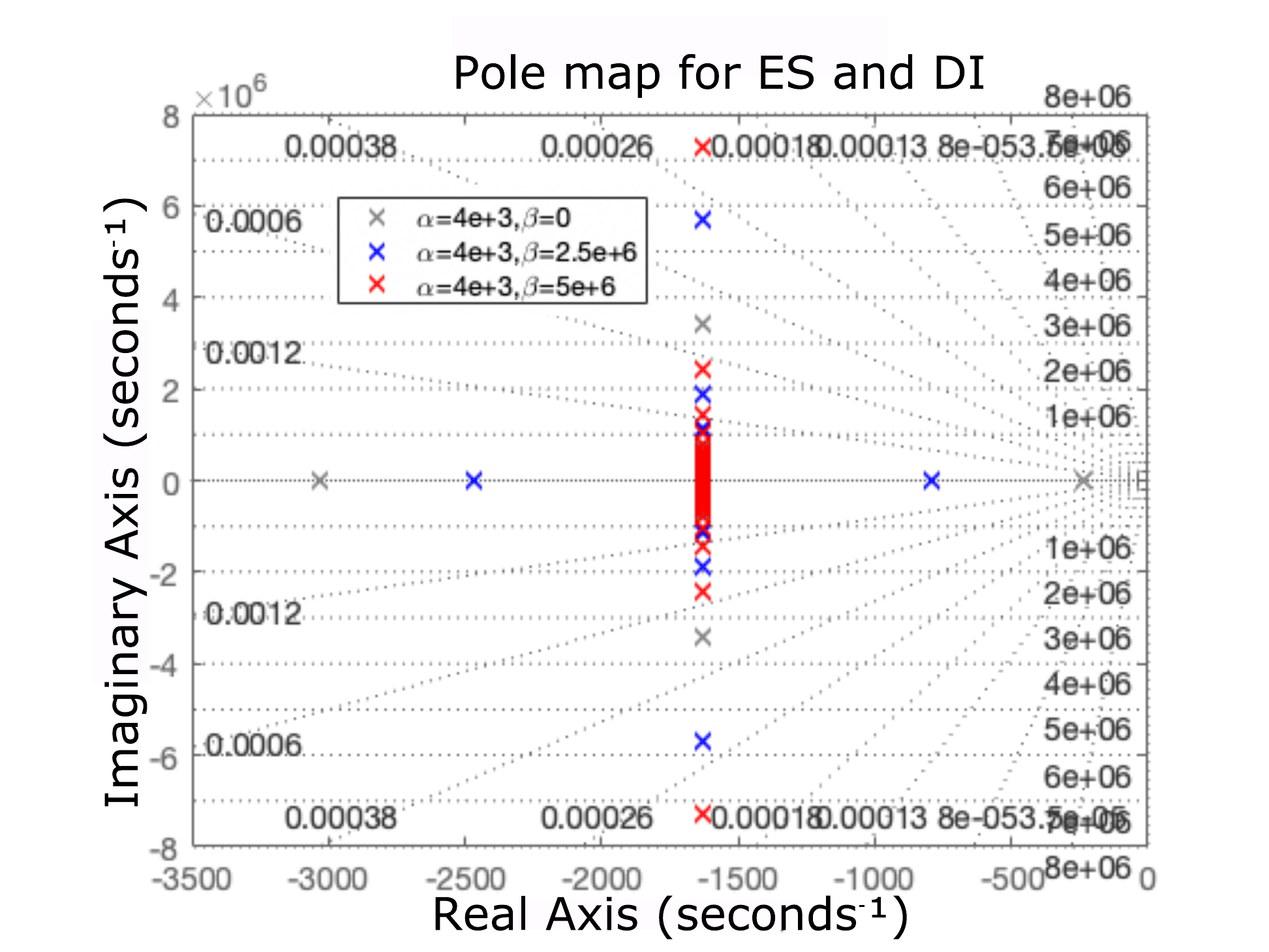}
			\caption{Poles of the closed-loop system with energy shaping and damping injection.}
			\label{fig:PolesZerosES-DI}
		\end{figure}

	Similarly, the poles of the closed-loop system with both the energy shaping and damping injection corresponding to Fig. \ref{fullyCBIdamping_position}(b) are presented in Fig. \ref{fig:PolesZerosES-DI}. The energy shaping has changed the natural frequency of the system and eliminated the aforementioned over-damping. Meanwhile, the bigger $\beta$ is, the higher frequency modes appear. But the stability of the closed-loop system is always guaranteed with the damping injection.
	
\begin{figure}[htbp!]
			\centering
		\subfloat[]{\includegraphics[scale=0.23]{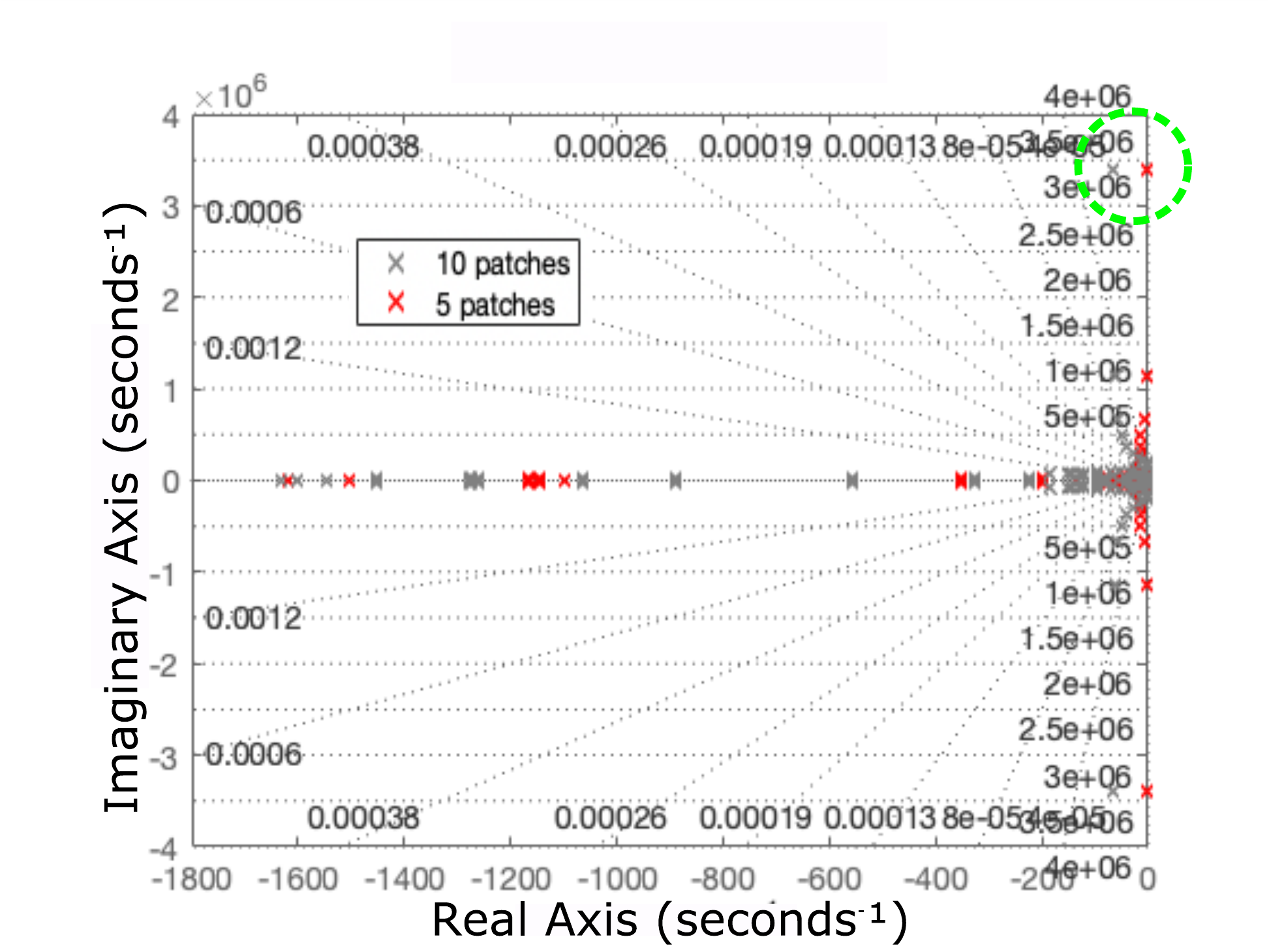}
		}
		\subfloat[]{\includegraphics[scale=0.23]{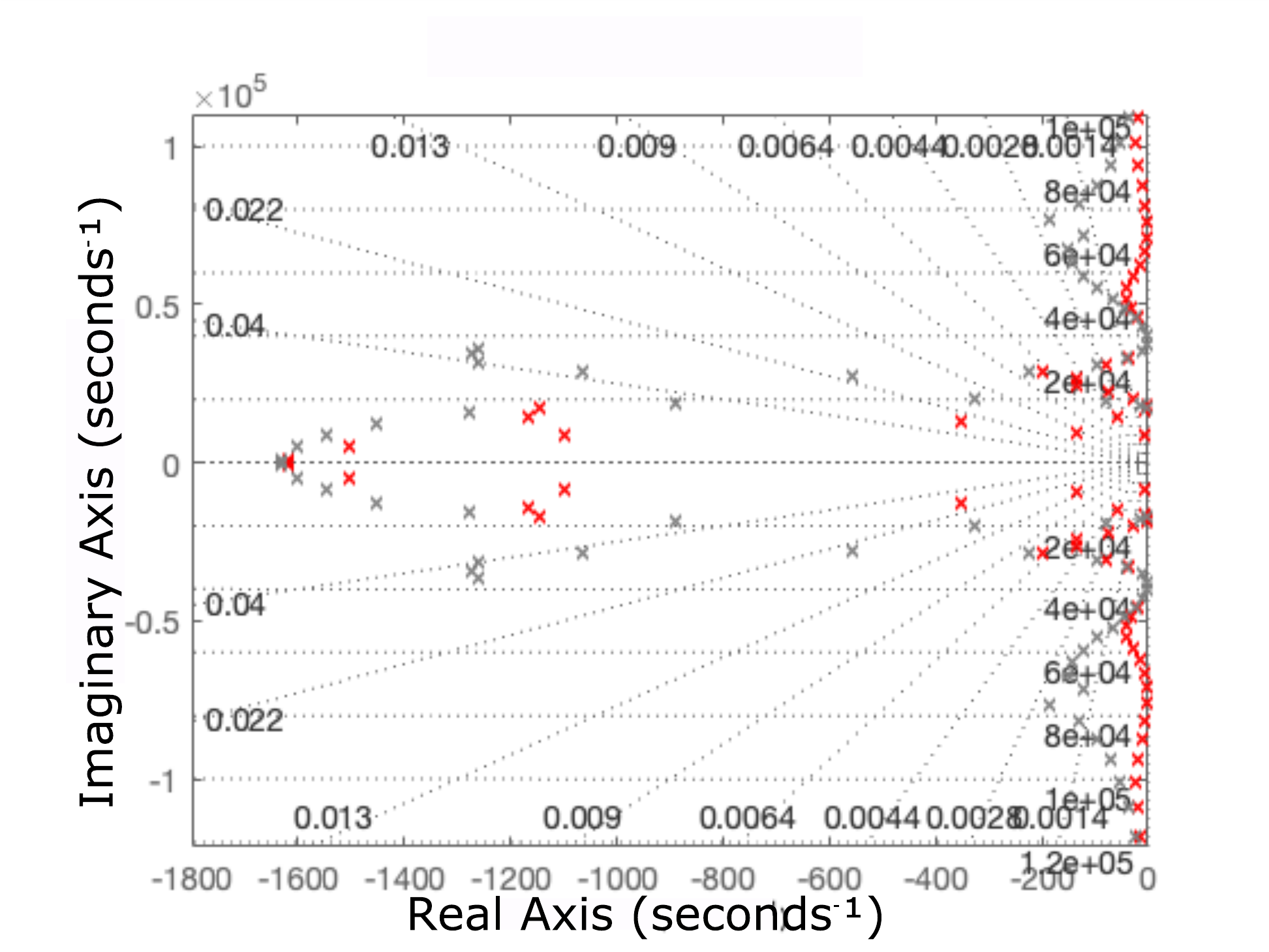}
		}
		\caption{Poles map for under-actuated case.
		}
		\label{fig:poleszerosunder}
	\end{figure}

	For the under-actuated case as presented in Fig. \ref{underactuateSim}, the poles are ploted in Fig. \ref{fig:poleszerosunder} (a), with its zoom of low frequency modes in Fig. \ref{fig:poleszerosunder} (b). It is shown that one can only control low frequency modes that dominate the response in order to approximate the poles distribution as in Fig. \ref{fig:PolesZerosES-DI} for fully-actuated case. The stability of the high frequency modes, e.g. the pole $s=-0.4\pm 3.4\times 10^6i$ of the closed-loop system with $5$ patches is preserved with the internal dissipation of the string.

\begin{figure}[htbp!]
	\centering
	\includegraphics[scale=.4]{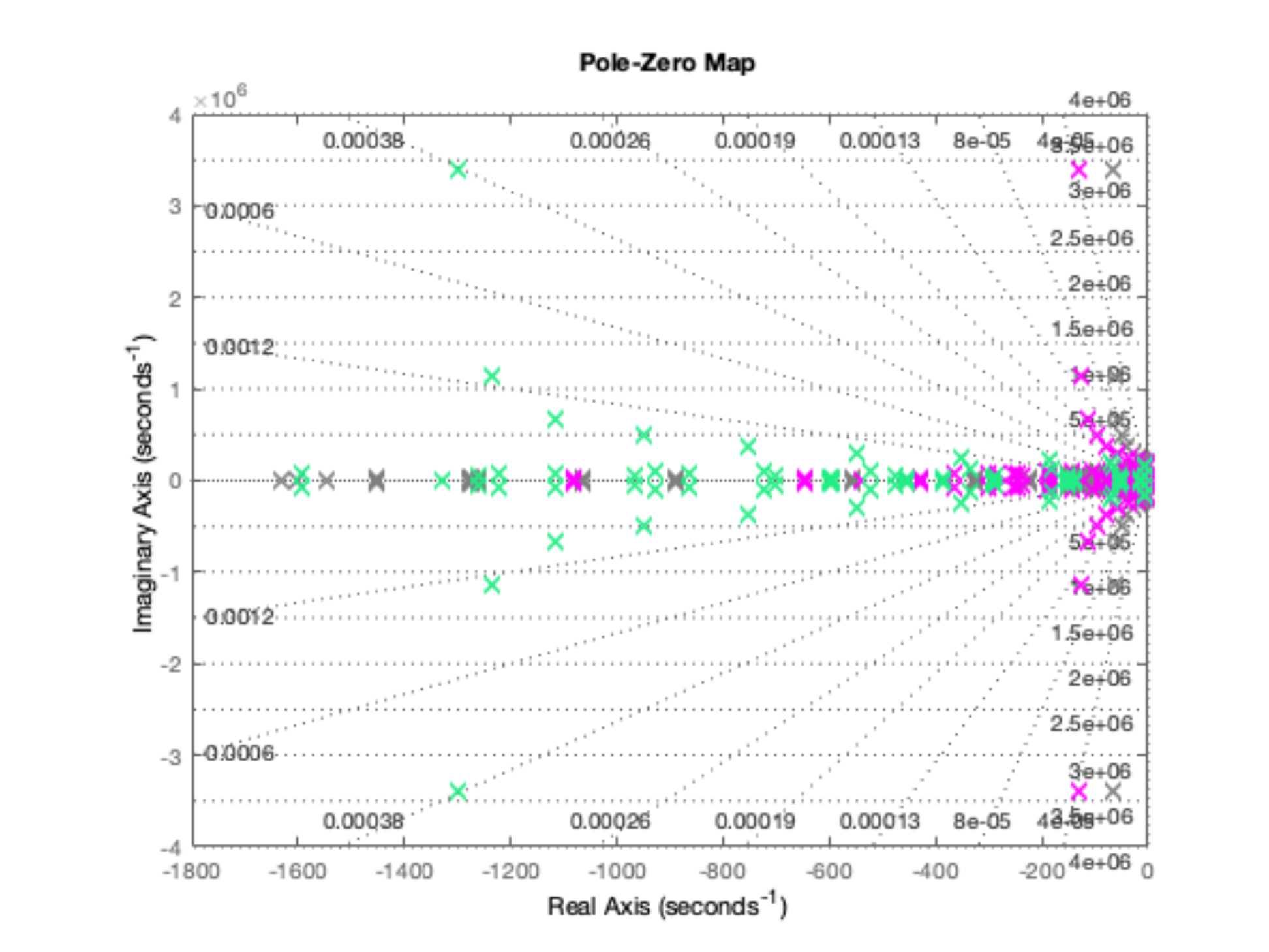}
	\caption{Poles of the closed-loop system for under-actuated case with the change of $D_c$.}
	\label{fig:PolesZerosunderChangeDc}
\end{figure}

\begin{remark}
	In order to investigate the influence of damping injection to the high frequency modes in under-actuated case, we vary the matrix $D_c$ and plot the closed-loop poles. According to Fig. \ref{fig:PolesZerosunderChangeDc},the conclusion is similar as the fully-actuated case presented in Fig.\ref{fig:PolesZerosDI}.
\end{remark}

\begin{figure}[htbp!]
	\centering
	\subfloat[]{\includegraphics[scale=0.23]{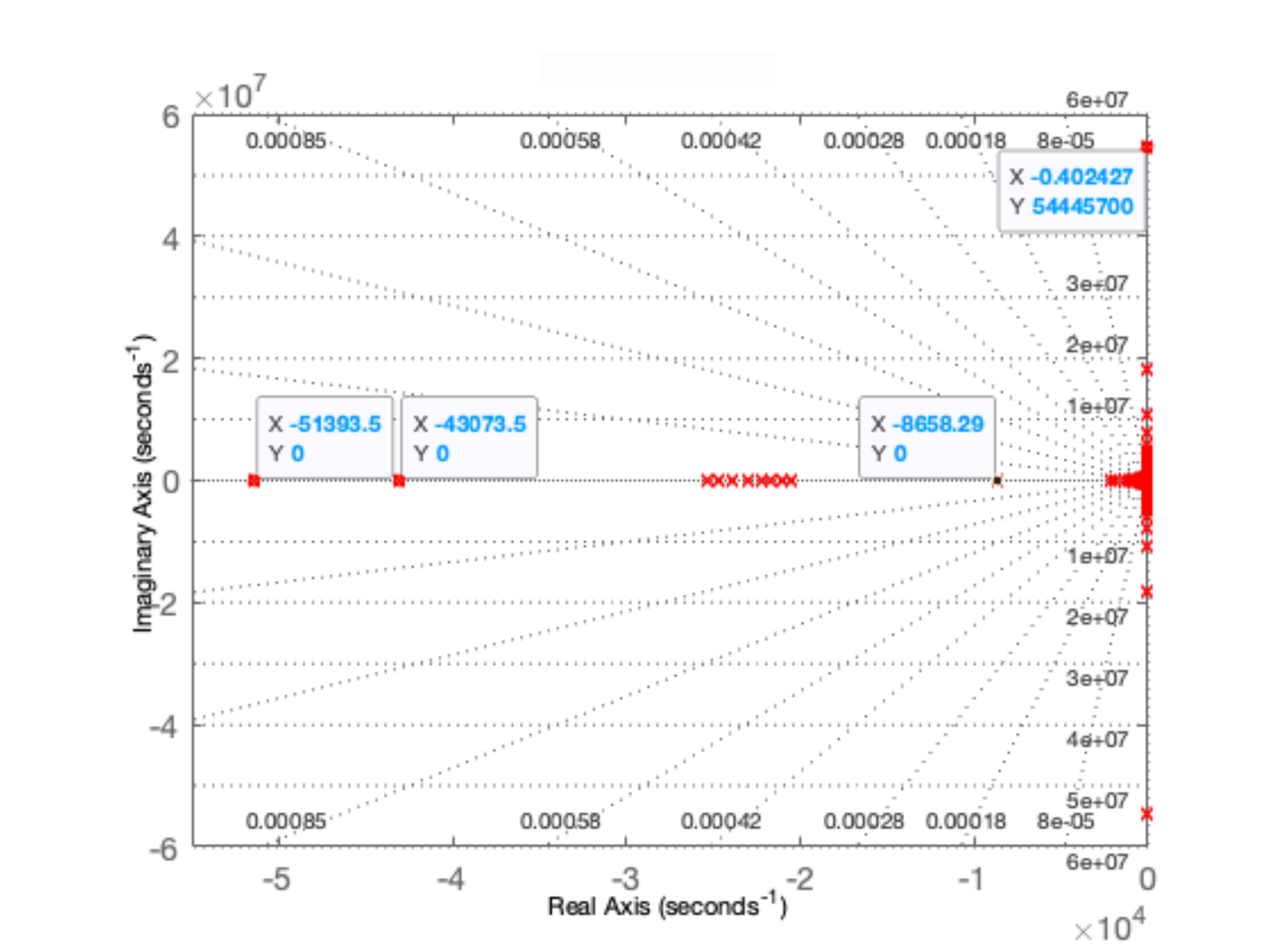}
	}
	\subfloat[]{\includegraphics[scale=0.23]{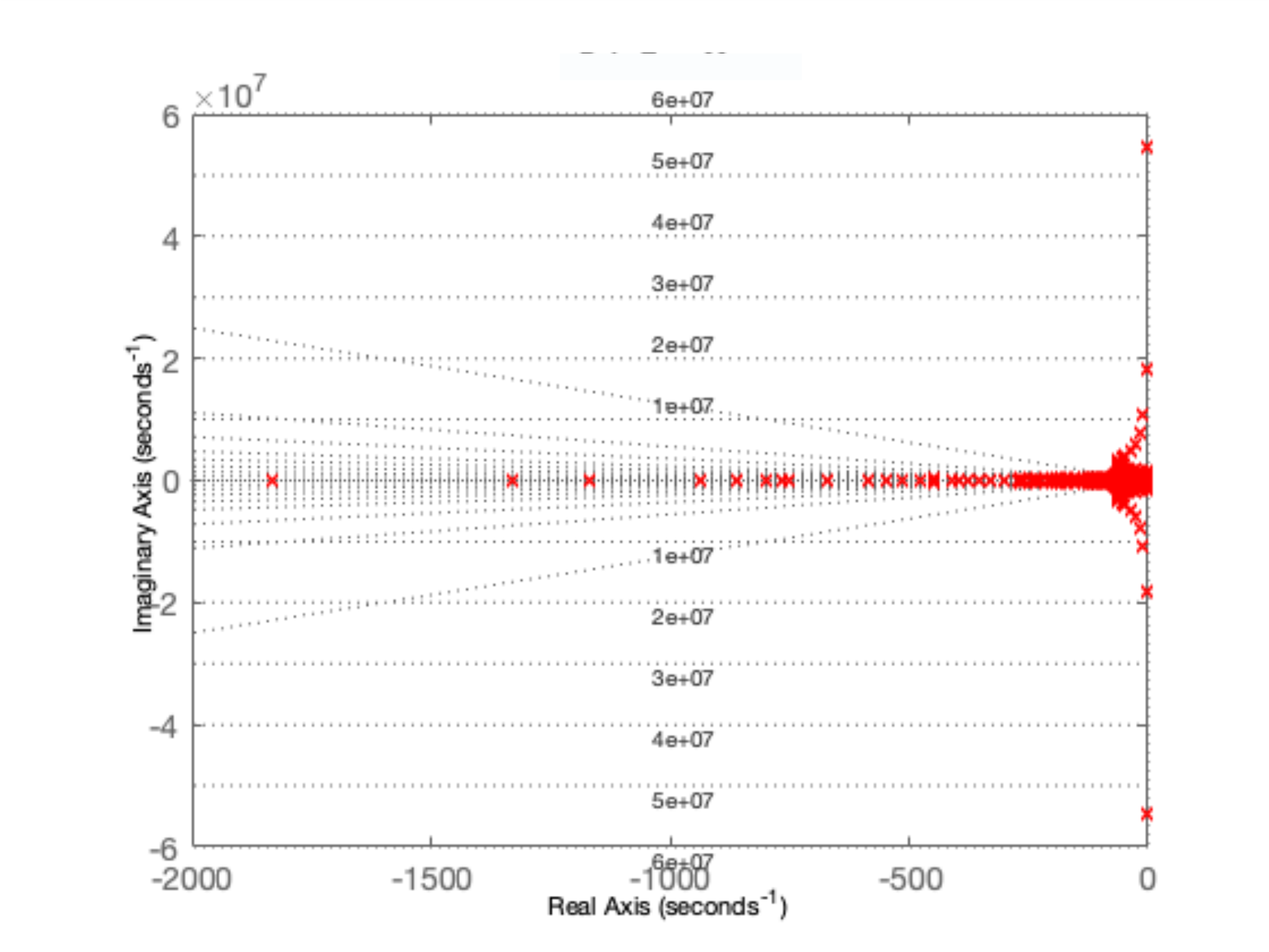}
	}
	\caption{Poles of the high-order closed-loop system with low-order designed controller.}
\label{fig:PoleZeroLow2High}
\end{figure}

When the designed low order controller with 10 patches is applied to a higher order system with $200$ discretized elements, the poles are presented in Figure \ref{fig:PoleZeroLow2High}. The stability of high frequency modes is always guaranteed with the internal dissipation. If we compare the poles at low frequency modes Fig. \ref{fig:PoleZeroLow2High}(a) with the Fig. \ref{fig:PolesZerosES-DI} and Fig. \ref{fig:poleszerosunder}(a), we can notice that these poles for high order system are more away from the imaginary axis, some of which are real poles $s=-51393.5,-43073.5$ and $-8658.29$. These poles have a very fast response such that they can be neglected. The zoom of high order poles in Fig. \ref{fig:PoleZeroLow2High}(b) are similar as in Fig. \ref{fig:poleszerosunder}(a) with $10$ actuator patches, which is consistent with the simulation results in Fig. \ref{ValidPrecise}(b).
	\section{Conclusion and future work}
	\label{SecConclusion}
	In this paper, we consider the in-domain control of infinite-dimensional port Hamiltonian systems with two conservation law using an early lumping approach. For control design purposes, we extend the CbI method to the use of controllers distributed in space. The distributed structural invariants are used to modify part of the closed-loop energy of the system. Two different cases are investigated: the ideal case where the system is fully-actuated and the under-actuated case where the control action is achieved using piecewise homogeneous inputs. In the latter the controller is derived by optimization. Simulations of both fully-actuated and under-actuated cases show how the damping injection together with the energy shaping improves the dynamic performances of the closed-loop system and keeps the closed-loop system asymptotically stable. Comparisons of the two cases also indicate that with an appropriate choice of the controller parameters, one can achieve similar performances for the under- and fully- actuated cases. Future works aim at extending the approach to the use of observers and at generalizing the proposed control design procedure to classes of non linear infinite-dimensional PHS.
	\bibliographystyle{plain}        
	\bibliography{biblio} 
\end{document}